\documentclass{amsart}
\usepackage{mathrsfs,yfonts}
\usepackage{amssymb,array,pb-diagram}

\newcounter{lemma}
\newtheorem{Theorem}{Theorem}
\newtheorem{theorem}{Theorem}

\newtheorem{Lemma}[lemma]{Lemma}
\newtheorem{Corollary}[lemma]{Corollary}
\newtheorem{Proposition}[lemma]{Proposition}
\newtheorem*{Example}{Example}
\newtheorem*{Remark}{Remark}
\newtheorem*{Conjecture}{Conjecture}

\def\F{\mathbb F}
\def\H{\mathbb H}

\def\Q{\mathbb Q}
\def\Z{\mathbb Z}

\def\C{\mathbb C}

\def\CC{\mathscr C}
\def\FF{\mathscr F}
\def\DD{\mathscr D}

\def\PP{\mathscr P}

\def\div{\operatorname{div}}
\def\Mod{\ \mathrm{mod}\ }
\def\gen#1{\langle#1\rangle}

\begin{document}

\title[Cuspidal rational torsion subgroup of $J_1(p^n)$]
  {Structure of the cuspidal rational torsion subgroup of $J_1(p^n)$}
\author{Yifan Yang}
\address{Department of Applied Mathematics \\
  National Chiao Tung University \\
  Hsinchu 300 \\
  TAIWAN}
\email{yfyang@math.nctu.edu.tw}
\author{Jeng-Daw Yu}
\address{Department of Mathematics and Statistics \\
  Queen's University \\
  Kingston, Canada K7L 3N6}
\email{jdyu@mast.queensu.ca}
\date\today

\begin{abstract} In this article, we determine the structure of the
  $p$-primary subgroup of the cuspidal rational torsion subgroup of
  the Jacobian $J_1(p^n)$ of the modular curve $X_1(p^n)$ for a
  regular prime $p$.
\end{abstract}

\maketitle

%%%%%%%%%%%%%%%%%%%%%%%%%%%%%%%%%%%%%%%%%%%%%%%%%%%%%%%%%%%%
%
%  Introduction
%
%%%%%%%%%%%%%%%%%%%%%%%%%%%%%%%%%%%%%%%%%%%%%%%%%%%%%%%%%%%%
\begin{section}{Introduction and statements of results} Let $\Gamma$
  be a congruence subgroup of $SL(2,\Z)$. The modular curve
  $X(\Gamma)$ and its Jacobian variety $J(\Gamma)$ are very important
  objects in number theory. For instance, the problem of determining
  all possible structures of ($\Q$-)rational torsion subgroup of
  elliptic curves over $\Q$ is equivalent to that of determining
  whether the modular curves $X_1(N)$ have non-cuspidal rational
  points. Also, the celebrated theorem of Wiles and others shows that
  every elliptic curve over $\Q$ is a factor of the Jacobian $J_0(N)$.
  In the present article, we are concerned with the arithmetic aspect
  of the Jacobian variety $J_1(N)$ of the modular curve $X_1(N)$. In
  particular, we will study the structure of the ($\Q$-)rational
  torsion subgroup of $J_1(N)$.

  Recall that the modular curve $X_1(N)$ possesses a model
  over $\Q$ on which the cusp $\infty$ is a ($\Q$-)rational point.
  (See \cite[Chapter 6]{Shimura} for details.)
  Thus, if $P$ is another rational cusp, then the image of $P$ under
  the cuspidal embedding $i_\infty:~X_1(N)\to J_1(N)$ sending $P$ to
  the divisor class $[(P)-(\infty)]$ will be a rational point on
  $J_1(N)$. Moreover, according to a result of Manin and Drinfeld
  \cite{Manin}, the point $i_\infty(P)$ is of finite order. In other
  words, the rational torsion subgroup of $J_1(N)$ contains a subgroup
  generated by the image of rational cusps under $i_\infty$, which we
  will refer to as the \emph{cuspidal rational torsion subgroup} of
  $J_1(N)$. In general, it is believed that the cuspidal rational
  torsion subgroup should be the whole rational torsion subgroup. (For
  primes $p$, the conjecture was formally stated in \cite[Conjecture
  6.2.2]{Conrad-Edixhoven-Stein}. The conjecture was verified for a
  few cases in the same paper.) Note that for the case $J_0(p)$, the
  Jacobian of $X_0(p)$ of prime level $p$, Mazur
  \cite[Theorem 1]{Mazur} has already shown that all rational torsion
  points are generated by the divisor class $[(0)-(\infty)]$.

  On the aforementioned model of $X_1(N)$, all the cusps of type $k/N$
  with $(k,N)=1$ are rational over $\Q$. (See, for example,
  \cite[Theorem 1.3.1]{Stevens-book}.)
  Moreover, if the level $N$ is relatively prime to $6$, then these
  cusps are the only rational cusps. Since these cusps are
  precisely those lying over $\infty$ of $X_0(N)$, for convenience, we
  shall call them the \emph{$\infty$-cusps}. Now suppose that we are
  given a divisor $D$ of degree $0$ on $X_1(N)$ supported on the
  $\infty$-cusps. Then the order of the divisor class $[D]$ in
  $J_1(N)$ is simply the smallest positive integer $m$ such that $mD$
  is a principal divisor, that is, the divisor of a modular function
  on $X_1(N)$. Therefore, to determine the group structure of the
  cuspidal rational torsion subgroup of $J_1(N)$, it is vital to
  study the group of modular units on $X_1(N)$ having divisors
  supported on the $\infty$-cusps. (In literature, if a modular
  function $f$ on a congruence subgroup $\Gamma$ has a divisor
  supported on cusps, then $f$ is called a \emph{modular unit}.)

  In a series of papers \cite{Kubert-Lang-I, Kubert-Lang-II,
  Kubert-Lang-III, Kubert-Lang-IV, Kubert-Lang-classno}, Kubert and
  Lang studied the group of modular units on $X(N)$ and $X_1(N)$. For
  the curves $X_1(N)$, they \cite[Chapter 3]{Kubert-Lang-book} showed
  that all modular units on $X_1(N)$ with divisors supported on the
  $\infty$-cusps are products of a certain class of Siegel functions.
  (See Section \ref{subsection: Siegel functions} for details.)
%  Moreover, they proved that the group of such modular units modulo
%  scalars has rank $\phi(N)/2-1$. Since the number of $\infty$-cusps is
%  $\phi(N)/2$, this means that every divisor of degree $0$ supported
%  on $\infty$-cusps has finite order in the Jacobian, in accordance
%  with Drinfeld and Manin's result.
  Furthermore, in
  \cite{Kubert-Lang-classno} they also determined the order of the
  torsion subgroup of $J_1(p^n)$ generated by the $\infty$-cusps for
  the case $p$ is a prime greater than $3$. (The case $N=p$ was first
  obtained in \cite{Klimek}.) Then Yu \cite{Yu} gave a formula for all
  positive integers $N$. (Note that all the results mentioned above
  dealt with modular units with divisors supported on the cusps lying
  over $0$ of $X_0(N)$, instead of the $\infty$-cusps, but it is easy
  to translate the results using the Atkin-Lehner involution
  $\left(\begin{smallmatrix}0&-1\\N&0\end{smallmatrix}\right)$.)

  In a very recent paper \cite{Yang2}, we applied Yu's divisor class
  number formula to determine an explicit basis for the group of
  modular units on $X_1(N)$ with divisors supported on the
  $\infty$-cusps for any positive integer $N$. As applications, we
  used the basis to compute the group structure of the cuspidal
  rational torsion subgroup of $J_1(N)$. A remarkable discovery is
  that when $p$ is a regular prime, the structure of the $p$-primary
  subgroup of the cuspidal rational torsion subgroup of $J_1(p^n)$
  seems to follow a simple pattern. (Recall that an odd prime $p$ is
  said to be \emph{regular} if $p$ does not divide the numerators of
  the Bernoulli numbers $B_2,B_4,\ldots,B_{p-3}$, or equivalently, if
  $p$ does not divide the class number of the cyclotomic field
  $\Q(e^{2\pi i/p})$.)

  More precisely, let $p$ be a prime, $n$ be a positive integer, and
  $\CC_1^\infty(p^n)$ be the subgroup of $J_1(p^n)$ generated by the
  $\infty$-cusps. Consider the endomorphism
  $[p]:\CC_1^\infty(p^n)\to\CC_1^\infty(p^n)$ defined by
  multiplication by $p$. Define the \emph{$p$-rank} of
  $\CC_1^\infty(p^n)$ to be the integer $k$ such that the kernel of
  $[p]$ has $p^k$ elements.

\begin{Conjecture}[Yang \cite{Yang2}] 
  Assume that $p$ is a regular prime. Then the $p$-rank of
  $\CC_1^\infty(p^n)$ is
  $$
    \frac12(p-1)p^{n-2}-1
  $$
  for prime power $p^n\ge 8$ with $n\ge 2$. More precisely, the number
  of copies of $\Z/p^{2k}\Z$ in the primary decomposition of
  $\CC_1^\infty(p^n)$ is
  $$
    \begin{cases}
    \frac12(p-1)^2p^{n-k-2}-1, &\text{if }p=2\text{ and }k\le n-3, \\
    \frac12(p-1)^2p^{n-k-2}-1, &\text{if }p\ge3\text{ and }k\le n-2, \\
    \frac12(p-5), &\text{if }p\ge 5\text{ and }k=n-1, \\
    0, &\text{else}.
    \end{cases}
  $$
  and the number of copies of $\Z/p^{2k-1}\Z$ is
  $$
    \begin{cases}
    1, &\text{if }p=2\text{ and }k\le n-3, \\
    1, &\text{if }p=3\text{ and }k\le n-2, \\
    1, &\text{if }p\ge 5\text{ and }k\le n-1, \\
    0, &\text{else}.
    \end{cases}
  $$
\end{Conjecture}

\begin{Example} \upshape For the primes $p=2$, $3$, $5$, the above
  conjecture asserts that the $p$-parts of $\CC_1^\infty(p^n)$ follow
  the pattern depicted in Table \ref{table 1}.
  \begin{table} \label{table 1}
\caption{$p$-primary part of $\CC_1^\infty(p^n)$}
$$ \extrarowheight3pt
\begin{array}{c||l} \hline\hline
p^n & p\text{-primary subgroups} \\ \hline\hline
2^4 & (2) \\ \hline
2^5 & (2)(2^2)(2^3) \\ \hline
2^6 & (2)(2^2)^3(2^3)(2^4)(2^5) \\ \hline
2^7 & (2)(2^2)^7(2^3)(2^4)^3(2^5)(2^6)(2^7)\\ \hline\hline
%2^8 & (2)(2^2)^{15}(2^3)(2^4)^7(2^5)(2^6)^3(2^7)(2^8)(2^9) \\ \hline\hline
3^3 & (3)(3^2) \\ \hline
3^4 & (3)(3^2)^5(3^3)(3^4)\\ \hline
3^5 & (3)(3^2)^{17}(3^3)(3^4)^5(3^5)(3^6) \\ \hline
3^6 & (3)(3^2)^{53}(3^3)(3^4)^{17}(3^5)(3^6)^5(3^7)(3^8)\\ \hline\hline
5^2 & (5) \\ \hline
5^3 & (5)(5^2)^7(5^3)\\ \hline
5^4 & (5)(5^2)^{39}(5^3)(5^4)^7(5^5)\\ \hline
5^5 & (5)(5^2)^{199}(5^3)(5^4)^{39}(5^5)(5^6)^7(5^7) \\ \hline\hline
\end{array}
$$
\end{table}
Here the notation $(p^{e_1})^{n_1}\ldots(p^{e_k})^{n_k}$ means that
the primary decomposition of $\CC_1^\infty(p^n)$ contains $n_i$ copies
of $\Z/p^{e_i}\Z$.
\end{Example}

The main purpose of the present article is to prove this conjecture.

\begin{Theorem} \label{theorem 1} The conjecture is true.
\end{Theorem}

We note that the assumption that $p$ is a regular prime is crucial in
the proof of Theorem \ref{theorem 1}. This assumption is used to
establish an exact formula for the $p$-rank of $\CC_1^\infty(p^n)$ and
to determine the kernel of the homomorphism
$\CC_1^\infty(p^n)\to\CC_1^\infty(p^{n-1})$ induced from the covering
$X_1(p^n)\to X_1(p^{n-1})$. At present, we do not know how to extend
our method to the case of irregular primes.

On the other hand, it is possible to obtain a similar result for
modular curves $X_1(p^nq^m)$, where $q$ is another prime, under the
assumption that the product
\begin{equation*}
  p\prod_{\chi}\frac14B_{2,\chi}
\end{equation*}
of generalized Bernoulli numbers $B_{2,\chi}$ associated with even
Dirichlet characters $\chi$ modulo $pq^m$ is a $p$-unit. For example,
following the argument in the present paper, we can show that the
$2$-primary subgroup of the torsion subgroup of $J_1(3\cdot 2^n)$
generated by the $\infty$-cusps is isomorphic to
$$
  \prod_{k=1}^{n-2}(\Z/2^{2k}\Z)^{2^{n-k-2}}.
$$
However, we will not pursue in this direction here because it does not
constitute a significant extension of Theorem \ref{theorem 1} and the
proof of some key lemmas in these cases is much more complicated
than the prime power cases. (For instance, it takes more than one page
just to describe the basis for the group of modular units on
$X_1(p^nq^m)$.)

The rest of the article is organized as follows. In Section
\ref{section: outline}, we describe our strategy in proving Theorem
\ref{theorem 1}. We will show that Theorem \ref{theorem 1} will follow
immediately from five properties of the divisor groups, namely,
Propositions \ref{proposition: C0}--\ref{proposition: kernel of p^2}.
In Section \ref{section: modular units}, we review our basis for the
group of modular units on $X_1(N)$, which constitutes the cornerstone
of our argument. In Section \ref{section: pi and iota},
we study the natural maps between the cuspidal groups
in different levels. We then give the proof of the five propositions
in Section \ref{section: propositions}.
\end{section}

%%%%%%%%%%%%%%%%%%%%%%%%%%%%%%%%%%%%%%%%%%%%%%%%%%%%%%%%%%%%
%
%  Outline
%
%%%%%%%%%%%%%%%%%%%%%%%%%%%%%%%%%%%%%%%%%%%%%%%%%%%%%%%%%%%%

\begin{section}{Outline of proof of Theorem \ref{theorem 1}}
  \label{section: outline} In this section, we will first collect all
  the notations and conventions used throughout the paper. We then
  describe our strategy in proving Theorem \ref{theorem 1}. Our
  arguments depend crucially on our explicit knowledge about the basis
  for the group of modular units on $X_1(N)$, which will be reviewed
  in Section \ref{subsection: basis}.

\begin{subsection}{Notations and conventions} \label{subsection: notations}
Let $p$ be a prime. We fix an integer $\alpha$ that generates
$(\Z/p^{n+1}\Z)^\times/\pm 1$ for all integers $n\ge 0$. Explicitly,
for $p=2$, we choose $\alpha=3$, and for an odd prime $p$,
we let $\alpha$ be an
integer that generates $(\Z/p\Z)^\times$, but satisfies
$\alpha^{p-1}\not\equiv 1\mod p^2$. For $n\ge 0$, we define
  \begin{align*}
    X_n&=\text{the modular curve }X_1(p^{n+1}), \\
    C_n&=\text{the set of cusps of }X_n
      \text{ lying over }\infty\text{ of }X_0(p^{n+1}),
      \text{ i.e., the set of }\infty\text{-cusps}, \\
    \phi_n&=|C_n|=\phi(p^{n+1})/2=p^n(p-1)/2, \\
    P_{n,k}&=\text{the cusp }\alpha^k/p^{n+1}\text{ in }C_n, \\
    \DD_n&=\text{the group of divisors of degree }0
      \text{ on }X_n\text{ having support on }C_n, \\
    \FF_n&=\text{the group of modular units on }X_n
      \text{ having divisors supported on }C_n, \\
    \PP_n&=\div\FF_n,\text{ the subgroup of principal divisors on }
      X_n\text{ having support on }C_n, \\
    \CC_n&=\DD_n/\PP_n, \text{ the rational torsion subgroup of
      }J_1(p^{n+1}) \text{ generated by }C_n, \\
    \pi_n&=\text{the canonical homomorphism from }\DD_n\text{ to }
      \DD_{n-1}\text{ induced from the covering } \\
    &\qquad\qquad X_n\to X_{n-1}. \\
    \iota_n&=\text{the embedding }\DD_{n-1}\to\DD_n\text{ defined by }
     \iota_n(P)=p\sum_{Q:~\pi_n(Q)=P}Q.
%    h_1^\infty(N)&=|\DD_1^\infty(N)/\div\FF_1^\infty(N)|,
%      \text{ the divisor class number}.
  \end{align*}
Note that $P_{n,k}$ and $P_{n,m}$ represent the same cusp on $X_n$ if
and only if $k\equiv m\mod\phi_n$. Then we have
$C_n=\{P_{n,k}:~k=0,\ldots,\phi_n-1\}$, and
$$
  \pi_n(P_{n,k})=P_{n-1,k}, \qquad
  \iota_n(P_{n-1,k})=p\sum_{h=0}^{p-1}P_{n,k+h\phi_{n-1}}.
$$

Since we are mainly interested in the orders of a function at the
$\infty$-cusps, for a modular function $f$ on $X_n$, we introduce the
notation $\div^\infty$ denoting the $C_n$-part
$$
  \div^\infty f=\sum_{P\in C_n}\operatorname{ord}_f(P)P
$$
of the divisor of $f$.

Finally, the generalized Bernoulli number $B_{k,\chi}$ associated
with a Dirichlet character $\chi$ modulo $N$, not necessarily
primitive, are defined by the series
$$
  \sum_{r=1}^N\frac{\chi(r)te^{rt}}{e^{Nt}-1}
 =\sum_{k=0}^\infty B_{k,\chi}\frac{t^k}{k!}.
$$
In particular, we have
$$
  B_{2,\chi}=N\sum_{r=1}^N\chi(r)B_2\left(\frac rN\right)
 =N\sum_{r=1}^N\chi(r)\left(\frac{r^2}{N^2}-\frac rN+\frac16\right).
$$
Here $B_2(x)=\{x\}^2-\{x\}+1/6$ and $\{x\}$ denotes the fractional
part of a real number $x$. The readers should be mindful that our
definition differs from some other authors' definition. See the remark
following Theorem \ref{theorem: Yu} for details.
\end{subsection}

\begin{subsection}{Outline of proof of Theorem \ref{theorem 1}}
\label{subsection: outline}
In this section, we will describe our strategy in proving Theorem
\ref{theorem 1}.

Intuitively, just by looking at Table \ref{table 1},
one immediately realizes that if the conjecture is to hold, then the
$p$-primary subgroup of $\CC_n/\ker[p^2]$ must have the same structure
as that of $\CC_{n-1}$, where $[p^2]$ denotes the
multiplication-by-$p^2$ homomorphism for an additive group, and one
expects that there should be a canonical isomorphism between the
$p$-primary subgroups of the two groups. The only sensible candidate
for such an isomorphism is the one induced by the covering
$X_n\to X_{n-1}$. To establish this isomorphism, we first show that
$\pi_n$ induces an isomorphism between the $p$-part of
$\CC_n=\DD_n/\PP_n$ and that of
$\pi_n(\DD_n)/\pi_n(\PP_n)=\DD_{n-1}/\pi_n(\PP_n)$. We then show
that the kernel of $[p^2]$ of the latter group is
$\PP_{n-1}/\pi_n(\PP_n)$, and thereby establish the isomorphism. The
following diagram illustrate the relations between various groups.
$$
  \begin{diagram}
  \node{\CC_n=\DD_n/\PP_n}\arrow{e,tb}{\pi_n}{p\text{-part }\simeq}
  \arrow{s,l}{/\ker[p^2]} \node{\DD_{n-1}/\pi_n(\PP_n)}
  \arrow{s,r}{/\ker[p^2]} \\
  \node{\CC_n/\ker[p^2]}\arrow{e,b}{p\text{-part }\simeq}
  \node{\CC_{n-1}=\DD_{n-1}/\PP_{n-1}}
  \end{diagram}
$$

Now assume that the isomorphism between the $p$-parts of
$\CC_n/\ker[p^2]$ and $\CC_{n-1}$ is established. This would show that
if the $p$-part of $\CC_{n-1}$ is $\prod(\Z/p^{e_i}\Z)^{r_i}$, then
the $p$-part of $\CC_n$ is
$(\Z/p\Z)^{s_1}\times(\Z/p^2\Z)^{s_2}\times\prod(\Z/p^{e_i+2}\Z)^{r_i}$
for some non-negative integers $s_1$ and $s_2$. If we can determine
the $p$-ranks of $\CC_{n-1}$ and $\CC_n$ and the index of
$\pi_n(\PP_n)$ in $\PP_{n-1}$, this will yield information about
$s_1+s_2$ and $s_1+2s_2$, respectively, which in turn will give us
the exact values of $s_1$ and $s_2$. Finally, if we know the structure
of $\CC_0$ ($\CC_1$ for $p=3$ and $\CC_2$ for $p=2$), then the
structure of the $p$-primary subgroup of $\CC_n$ is determined for all
$n$.

In summary, to establish Theorem \ref{theorem 1}, it suffices to prove
the following propositions.

\begin{Proposition} \label{proposition: C0} If $p$ is a
  regular prime, then $p$ does not divide $|\CC_0|$. Also, if $p=2,3$,
  then $p\nmid|\CC_1|$, and if $p=2$, then $p\nmid|\CC_2|$.
\end{Proposition}

\begin{Proposition} \label{proposition: p-rank} Let $p$ be a regular
  prime. If $p^{n+1}\ge 5$, then the $p$-rank of $\CC_n$ is
  $p^{n-1}(p-1)/2-1$.
\end{Proposition}

\begin{Proposition} \label{proposition: index of projection} For all
  primes $p$, we have $\pi_n(\PP_n)\subset\PP_{n-1}$,
  and the index of $\pi_n(\PP_n)$ in $\PP_{n-1}$ is
  $p^{p^{n-1}(p-1)-3}$ if $p^{n+1}\ge 5$. Moreover, the structure of
  the factor group $\PP_{n-1}/\pi_n(\PP_n)$ is
  $$
    (\Z/p^2\Z)^{p^{n-1}(p-1)/2-2}\times(\Z/p\Z).
  $$
\end{Proposition}

\begin{Proposition} \label{proposition: projection is an isomorphism}
  Assume that $p$ is a regular prime. Then the $p$-part of
  $\CC_n=\DD_n/\PP_n$ is isomorphic to the $p$-part of
  $\DD_{n-1}/\pi_n(\PP_n)$.
\end{Proposition}

\begin{Proposition} \label{proposition: kernel of p^2} Let $p$ be a
  prime. Then the kernel of the multiplication-by-$p^2$ endomorphism
  $[p^2]$ of $\DD_{n-1}/\pi_n(\PP_n)$ is $\PP_{n-1}/\pi_n(\PP_n)$.
\end{Proposition}

\begin{Remark} \upshape We remark that the assumption that $p$ is a
  regular prime is crucial in the proof of Propositions
  \ref{proposition: C0}, \ref{proposition: p-rank}, and
  \ref{proposition: projection is an isomorphism}. In fact, the
  assumption is a necessary and sufficient condition for the three
  propositions. For example, by carefully examining the proof of
  Proposition \ref{proposition: p-rank}, one sees that if $p$ is an
  irregular prime, then the $p$-rank of $\CC_n$ is strictly greater
  than $p^{n-1}(p-1)/2-1$.

  Note also that Propositions \ref{proposition: projection is an
  isomorphism} and \ref{proposition: kernel of p^2} together imply
  that when $p$ is a regular prime, the $p$-part of $\CC_n/\ker[p^2]$
  is isomorphic to that of $\CC_{n-1}$. In terms of Jacobians, the
  $p$-part of $\CC_n$ corresponds to the rational $p$-power-torsion
  subgroup of $J_1(p^{n+1})$. So what these two propositions really
  say is that when $p$ is a regular prime, the kernel of the canonical
  homomorphism $\pi:J_1(p^{n+1})\to J_1(p^n)$ agrees with the kernel
  of $[p^2]:J_1(p^{n+1})\to J_1(p^{n+1})$ on the cuspidal rational
  $p$-power-torsion part of $J_1(p^{n+1})$, that is, on the cuspidal
  part,
  $$
    (p\text{-power torsion})\cap\ker\pi
   =(p\text{-power torsion})\cap\ker[p^2].
  $$
  Note that $p^2$ is exactly the degree of the covering $X_n\to
  X_{n-1}$.
  Naturally, one wonders whether it is still the case when $p$ is an
  irregular prime. We do not know the answer to this question.
\end{Remark}

The proof of these propositions will be postponed until Section
\ref{section: propositions}. Here let us formally complete the
proof of Theorem \ref{theorem 1}, assuming the truth of the
propositions.

\begin{proof}[Proof of Theorem \ref{theorem 1}]
  By Propositions \ref{proposition: projection is an isomorphism} and
  \ref{proposition: kernel of p^2}, when $p$ is a regular prime,
  \begin{equation*}
  \begin{split}
    p\text{-part of }\CC_n/\ker[p^2] & \simeq
    p\text{-part of }(\DD_{n-1}/\pi_n(\PP_n))/\ker[p^2] \\
   &=p\text{-part of }(\DD_{n-1}/\pi_n(\PP_n))\Big/
    (\PP_{n-1}/\pi_n(\PP_n)) \\
   &\simeq p\text{-part of }\DD_{n-1}/\PP_{n-1}=\CC_{n-1},
  \end{split}
  \end{equation*}
  Thus, if the structure of the $p$-part of $\CC_{n-1}$ is
  $$
    \prod_{i=1}^k(\Z/p^{e_i}\Z)^{r_i},
  $$
  then according to the structure theorem for finite abelian groups,
  the structure of the $p$-part of $\CC_n$ is
  $$
    (\Z/p\Z)^{s_1}\times(\Z/p^2\Z)^{s_2}\times
    \prod_{i=1}^k(\Z/p^{e_i+2}\Z)^{r_i}
  $$
  for some non-negative integers $s_1$ and $s_2$. Here the sum of
  $r_i$ is what we call the $p$-rank of $\CC_{n-1}$, and the sum
  of $s_1$, $s_2$, and $r_i$ is the $p$-rank of $\CC_n$. Using
  Proposition \ref{proposition: p-rank}, we find the integers $s_1$
  and $s_2$ satisfy
  \begin{equation} \label{equation: proof of Theorem 1 1}
    s_1+s_2=\frac12p^{n-2}(p-1)^2.
  \end{equation}
  On the other hand, by Propositions \ref{proposition: index of
  projection} and \ref{proposition: projection is an isomorphism}, we
  know that
  \begin{equation*}
  \begin{split}
    p\text{-part of }|\CC_n|/|\CC_{n-1}|=|\PP_{n-1}/\pi_n(\PP_n)|
   =p^{p^{n-1}(p-1)-3},
  \end{split}
  \end{equation*}
  which, together with Proposition \ref{proposition: p-rank}, implies
  that
  $$
    s_1+2s_2=(p^{n-1}(p-1)-3)-2(p^{n-2}(p-1)/2-1)=p^{n-2}(p-1)^2-1.
  $$
  Combining this with \eqref{equation: proof of Theorem 1 1}, we get
  $s_1=1$ and $s_2=p^{n-2}(p-1)^2/2-1$. Finally, Proposition
  \ref{proposition: C0} shows that the $p$-part of $\CC_0$ ($\CC_1$
  for $p=3$ and $\CC_2$ for $p=2$) is trivial. Then an induction
  argument gives the claimed result.
\end{proof}
\end{subsection}
\end{section}

%%%%%%%%%%%%%%%%%%%%%%%%%%%%%%%%%%%%%%%%%%%%%%%%%%%%%%%%%%%%
%
%  Modular units
%
%%%%%%%%%%%%%%%%%%%%%%%%%%%%%%%%%%%%%%%%%%%%%%%%%%%%%%%%%%%%

\begin{section}{Group of modular units on $X_1(N)$}
\label{section: modular units}
In this section, we will introduce our basis for the group $\FF_n$,
which is essential in our proof of Theorem \ref{theorem 1}. The
construction of our basis utilizes the Siegel functions.

\begin{subsection}{Siegel functions}
\label{subsection: Siegel functions}
  The Siegel functions are usually defined as products of the Klein
  forms and the Dedekind eta function. For our purpose, we only need
  to know that they have the following infinite product
  representation.

  For a pair of rational numbers $(a_1,a_2)\in\Q^2\backslash\Z^2$ and
  $\tau\in\H$, set $z=a_1\tau+a_2$, $q_\tau=e^{2\pi i\tau}$, and
  $q_z=e^{2\pi iz}$.
  Then the Siegel function $G_{(a_1,a_2)}(\tau)$ satisfies
  $$
    G_{(a_1,a_2)}(\tau)=-e^{2\pi ia_2(a_1-1)/2}q_\tau^{B(a_1)/2}
    (1-q_z)\prod_{n=1}^\infty(1-q_\tau^n q_z)(1-q_\tau^n/q_z),
  $$
  where $B(x)=x^2-x+1/6$ is the second Bernoulli polynomial. To
  construct modular units on $X_1(N)$ with divisors supported on the
  $\infty$-cusps, we consider a special class of Siegel functions.

  Given a positive integer $N$ and an integer $a$ not congruent to $0$
  modulo $N$, we set
  $$
    E_a^{(N)}(\tau)=-G_{(a/N,0)}(N\tau)
   =q^{NB(a/N)/2}\prod_{n=1}^\infty\left(1-q^{(n-1)N+a}\right)
    \left(1-q^{nN-a}\right),
  $$
  where $q=e^{2\pi i\tau}$. If the integer $N$ is clear from the
  context, we will write $E_a$ in place of $E_a^{(N)}$.

  We now review the properties of $E_a$. The material is mainly taken
  from \cite{Yang}. For more details, see op. cit. In the first lemma,
  we describe two simple, but yet very important relations between
  Siegel functions of two different levels.

  \begin{Lemma} \label{lemma: Eg relations} Let $M$ and $N$ be two
  positive integer. Assume that $N=nM$ for some integer $n$.
  Let $a$ be an
  integer not congruent to $0$ modulo $N$. Then
  \begin{equation} \label{equation: EN EM}
    E^{(N)}_{na}(\tau)=E^{(M)}_a(n\tau).
  \end{equation}
  Moreover, we have for all integers $a$ with $0<a<M$,
  \begin{equation} \label{equation: Bernoulli relation}
    N\sum_{k=0}^{n-1}B_2\left(\frac{kM+a}N\right)
   =MB_2\left(\frac aM\right),
  \end{equation}
  and consequently
  \begin{equation} \label{equation: EN EM 2}
    \prod_{k=0}^{n-1}E_{kM+a}^{(N)}(\tau)=E_a^{(M)}(\tau).
  \end{equation}
  \end{Lemma}

  \begin{proof} Relation \eqref{equation: EN EM} follows trivially
    from the definition of $E^{(N)}_g$. Property \eqref{equation:
    Bernoulli relation} can be verified by a direct computation.
    Relation \eqref{equation: EN EM 2} is an immediate consequence of
    \eqref{equation: Bernoulli relation} and the definition of
    $E^{(N)}_a$.
  \end{proof}

  The next lemma gives the transformation law for $E_a$ under the
  action of matrices in $\Gamma_0(N)$.

\begin{Lemma}[{\cite[Corollary 2]{Yang}}]
  \label{lemma: Eg formulas} For integers $g$ not congruent to $0$
  modulo $N$, the functions $E_g$ satisfy
\begin{equation}
\label{shifting for Eg}
  E_{g+N}=E_{-g}=-E_g.
\end{equation}
Moreover, let $\gamma=\begin{pmatrix}a&b\\ cN&d\end{pmatrix}\in
\Gamma_0(N)$. We have, for $c=0$,
$$
  E_g(\tau+b)=e^{\pi ibNB(g/N)}E_g(\tau),
$$
and, for $c>0$,
\begin{equation}\label{Eg formula}
  E_g(\gamma\tau)=\epsilon(a,bN,c,d)e^{\pi i(g^2ab/N-gb)}
  E_{ag}(\tau),
\end{equation}
where
\begin{equation*}
  \epsilon(a,b,c,d)
  =\begin{cases}
    e^{\pi i\left(bd(1-c^2)+c(a+d-3)\right)/6},
      &\text{if }c\text{ is odd}, \\
    -ie^{\pi i\left(ac(1-d^2)+d(b-c+3)\right)/6},
      &\text{if }d\text{ is odd}.
  \end{cases}
\end{equation*}
\end{Lemma}

\begin{Remark} \upshape
Note that Property \eqref{shifting for Eg} implies that there are only
$\lceil(N-1)/2\rceil$ essentially distinct $E_g$, indexed over the set
$(\Z/N\Z)/\pm 1-\{0\}$. Hence, a product $\prod_g$ or a sum $\sum_g$
is understood to be running over $g\in(\Z/N\Z)^\times/\pm 1$.
\end{Remark}

The functions $E_g$ clearly have no poles nor zeros in the upper
half-plane. The next lemma describes the order of $E_g$ at cusps of
$X_1(N)$.

\begin{Lemma}[{\cite[Lemma 2]{Yang}}]
\label{lemma: behavior of Eg}
The order of the function $E_g$ at a cusp $a/c$ of $X_1(N)$ with
$(a,c)=1$ is $(c,N)B_2(ag/(c,N))/2$, where $B_2(x)=\{x\}^2-\{x\}+1/6$
and $\{x\}$ denotes the fractional part of a real number $x$.
\end{Lemma}

The following theorem of Yu \cite{Yu} characterizes the modular units
on $X_1(N)$ with divisors supported at the $\infty$-cusps in terms of
$E_g$.

\begin{theorem}[{\cite[Theorem 4]{Yu}}] \label{theorem: DQ DO} Let $N$
  be a positive integer. A modular function $f$ on $\Gamma_1(N)$ has a
  divisor supported on the cusps $k/N$, $(k,N)=1$, if and only if
  $f=\prod_g E_g^{e_g}$ with the exponents $e_g$ satisfying the two
  conditions
  \begin{equation} \label{condition DQ}
    \sum_g g^2e_g\equiv 0\mod
    \begin{cases} N, &\text{if }N\text{ is odd}, \\
    2N, &\text{if }N\text{ is even}, \end{cases}
  \end{equation}
  and
  \begin{equation} \label{condition DO}
    \sum_{g\equiv\pm a\Mod N/p}e_g=0
  \end{equation}
  for all prime factors $p$ of $N$ and all integers $a$.
\end{theorem}

We remark that, again, Theorem 4 of \cite{Yu} was stated in the
setting of modular units with divisor supported on the $0$-cusps,
i.e., the cusps lying over $0$ of $X_0(N)$. Here we use the
Atkin-Lehner involution
$\left(\begin{smallmatrix}0&-1\\N&0\end{smallmatrix}\right)$ to
get Theorem \ref{theorem: DQ DO} from Yu's result.
\end{subsection}

\begin{subsection}{Basis for $\FF_n$} \label{subsection: basis}
  We now describe our basis for $\FF_n$ constructed in \cite{Yang2}.
  The case of an odd prime $p$ and the case of $p=2$ are stated in
  Theorems \ref{theorem: prime power} and \ref{theorem: power of 2},
  respectively.

\begin{theorem}[{\cite[Theorem 2]{Yang2}}] \label{theorem: prime power}
  Let $n\ge0$ and $N=p^{n+1}$ be an odd prime power. For a
  non-negative integer $\ell$, we set $\phi_\ell=\phi(p^{\ell+1})/2$.
  Let $\alpha$ be a generator of the cyclic group
  $(\Z/p^{n+1}\Z)^\times/\pm 1$, and $\beta$ be an integer such that
  $\alpha\beta\equiv 1\mod p$. Then a basis for $\FF_n$ modulo
  $\C^\times$ is given by
$$
  \begin{cases}
  f_i=\frac{E_{\alpha^{i-1}}E_{\alpha^{i+\phi_{n-1}}}^{\beta^2}}
  {E_{\alpha^{i+\phi_{n-1}-1}}E_{\alpha^i}^{\beta^2}},
  &i=1,\ldots,\phi_n-\phi_{n-1}-1, \\
  f_i=\frac{E_{\alpha^{i-1}}^p}{E_{\alpha^{i+\phi_{n-1}-1}}^p},
  &i=\phi_n-\phi_{n-1}, \\
  f_i=\frac{E^{(p^n)}_{\alpha^{i-1}}(p\tau)}
  {E^{(p^n)}_{\alpha^{i+\phi_{n-2}-1}}(p\tau)},
  &i=\phi_n-\phi_{n-1}+1,\ldots,\phi_n-\phi_{n-2}, \\
  \qquad \vdots & \qquad \vdots \\
  f_i=\frac{E^{(p^2)}_{\alpha^{i-1}}(p^{n-1}\tau)}
    {E^{(p^2)}_{\alpha^{i+\phi_0-1}}(p^{n-1}\tau)},
  &i=\phi_n-\phi_1+1,\ldots,\phi_n-\phi_0, \\
  f_i=\frac{E^{(p)}_{\alpha^{i-1}}(p^n\tau)}
  {E^{(p)}_{\alpha^i}(p^n\tau)}, &i=\phi_n-\phi_0+1,\ldots,\phi_n-1.
  \end{cases}
$$
\end{theorem}

\begin{theorem}[{\cite[Theorem 3]{Yang2}}] \label{theorem: power of 2}
  Let $n\ge 2$ and $N=2^{n+1}$. Let $\alpha=3$ be a generator of the
  cyclic group $(\Z/2^{n+1}\Z)^\times/\pm1$. 
  For $\ell\ge 1$, set $\phi_\ell=\phi(2^{\ell+1})/2=2^{\ell-1}$. Then
  a basis for $\FF_n$ modulo $\C^\times$ is given by
$$
  \begin{cases}
  f_i=\frac{E_{\alpha^{i-1}}E_{\alpha^{i+\phi_{n-1}}}}
    {E_{\alpha^i}E_{\alpha^{i+\phi_{n-1}-1}}},
    &i=1,\ldots,\phi_n-\phi_{n-1}-1, \\
  f_i=\frac{E_{\alpha^{i-1}}^2}{E_{\alpha^{i+\phi_{k-1}-1}}^2},
    &i=\phi_n-\phi_{n-1}, \\
  f_i=\frac{E_{\alpha^{i-1}}^{(2^n)}(2\tau)}
    {E_{\alpha^{i+\phi_{n-2}-1}}^{(2^n)}(2\tau)},
    &i=\phi_n-\phi_{n-1}+1,\ldots,\phi_n-\phi_{n-2}, \\
  \qquad\vdots & \qquad\vdots \\
  f_i=\frac{E_{\alpha^{i-1}}^{(8)}(2^{n-2}\tau)}
    {E_{\alpha^i}^{(8)}(2^{n-2}\tau)}, &i=\phi_n-1.
  \end{cases}
$$
\end{theorem}

The proof of these two theorems use the following divisor class number
formula of Kubert, Lang, and Yu, which will also be used in the
present paper. Note that the cases $p\ge 5$ were proved in
\cite{Kubert-Lang-classno}, while the cases $p=2,3$ were settled in
\cite{Yu}. In the same paper \cite{Yu}, Yu also obtained a divisor
class number formula for general $N$, although the general result is
not needed in the present article.

\begin{theorem}[{\cite[Theorem 3.4]{Kubert-Lang-classno} and
    \cite[Theorem 5]{Yu}}] \label{theorem: Yu}
  Let $N=p^{n+1}$ be a prime power greater than $4$. We have the
  divisor class number formula
  \begin{equation} \label{equation: Yu's formula}
    |\CC_n|=p^{L(p)}
    \prod_{\chi\neq\chi_0\text{ even}}\frac14B_{2,\chi},
  \end{equation}
  where
  $$
    L(p)=\begin{cases}
    p^{n-1}-2n+2, &\text{if }N=p^n\text{ and }p\text{ is odd}, \\
    2^{n-1}-2n+3, &\text{if }N=2^n\ge 8,
    \end{cases}
  $$
  and the product runs over all even non-principal
  Dirichlet characters modulo $p^{n+1}$.
\end{theorem}

\begin{Remark} \upshape We should remark that the definition of
  generalized Bernoulli numbers used in \cite{Kubert-Lang-classno} and
  \cite{Yu} is different from ours. Namely, if an even Dirichlet
  character $\chi$ modulo $N$ has a conductor $f$, then their
  definition is given by
  $$
    \frac12\sum_{r=1}^f\frac{\chi_f(r)te^{rt}}{e^{ft}-1}
   =\sum_{k=0}^\infty B_{2,\chi}\frac{t^k}{k!},
  $$
  where $\chi_f$ is the Dirichlet character modulo $f$ that induces
  $\chi$. When $N$ is a prime power $p^n$ and $\chi$ is not principal,
  the two definitions differ by a $1/2$ factor.

  Moreover, the readers are reminded that there were slight errors in
  the original statement of \cite[Theorem 5]{Yu}. See the discussion
  following Theorem A of \cite{Yang2} for details.
\end{Remark}
\end{subsection}
\end{section}

%%%%%%%%%%%%%%%%%%%%%%%%%%%%%%%%%%%%%%%%%%%%%%%%%%%%%%%%%%%%
%
%  Properties of pi and iota
%
%%%%%%%%%%%%%%%%%%%%%%%%%%%%%%%%%%%%%%%%%%%%%%%%%%%%%%%%%%%%

\begin{section}{Properties of $\pi_n$ and $\iota_n$}
\label{section: pi and iota}
  Throughout the section, we will follow the notations specified in
  Section \ref{subsection: notations}. The main results in this
  section are Lemmas \ref{lemma: projection principal to principal}
  and \ref{lemma: iota 2}, which state that $\pi_n$ maps a principal
  divisor to a principal divisor, and that if $\iota_n(D)$ is a
  principal divisor, then $D$ itself is principal. In addition, in
  Lemma \ref{lemma: iota 1} we will prove the converse to Lemma
  \ref{lemma: iota 2}, that is, if $D$ is a principal divisor in
  $\DD_{n-1}$, then $\iota_n(D)$ is a principal divisor.
%This property
%  is not used in the proof of Theorem \ref{theorem 1}.
% However, Lemmas
%  \ref{lemma: iota 1} and \ref{lemma: iota 2} combined give an
%  embedding of $\CC_{n-1}$ into $\CC_n$.

  The first lemma is rather trivial, but it plays a crucial role in
  the proof of Proposition \ref{proposition: kernel of p^2}.

\begin{Lemma} \label{lemma: p^2} We have
\begin{equation} \label{equation: p^2}
  \pi_n\circ\iota_n=[p^2],
\end{equation}
the multiplication-by-$p^2$ endomorphism of $\DD_{n-1}$
\end{Lemma}

\begin{proof} Obvious.
\end{proof}

In the next lemma we compute the image of the divisor of
$E^{(p^{n+1})}_g$ under $\pi_n$. Here we recall that the notation
$\div^\infty f$ means the $C_n$-part of the divisor of $f$.

\begin{Lemma} \label{lemma: divisor of Eg} Let $g$ be an integer. For
  $g\not\equiv 0\mod p^{n+1}$, we have
$$
  \pi_n(\div^\infty E_g^{(p^{n+1})})=
  \begin{cases}\div^\infty E_g^{(p^n)}, &\text{if }p\nmid g, \\
  p^2\div^\infty E_{g/p}^{(p^n)}, &\text{if }p|g.\end{cases}
$$
\end{Lemma}

\begin{proof} By Lemma \ref{lemma: behavior of Eg}, we have
$$
  \div^\infty E_g^{(p^{n+1})}=\frac{p^{n+1}}2\sum_{k=0}^{\phi_n-1}
  B_2\left(\frac{g\alpha^k}{p^{n+1}}\right)P_{n,k}.
$$
Recall that $\pi_n(P_{n,k})=\pi_n(P_{n,h})$ if and only if $h\equiv
k\mod\phi_{n-1}$. Thus,
$$
  \pi_n(\div^\infty E_g^{(p^{n+1})})=\frac{p^{n+1}}2
  \sum_{k=0}^{\phi_{n-1}-1}P_{n-1,k}\sum_{h=0}^{p-1}
  B_2\left(\frac{g\alpha^{k+h\phi_{n-1}}}{p^{n+1}}\right).
$$
Now assume that $p$ does not divide $g$, then as $h$ goes through $0$
to $p-1$, the residue classes of $g\alpha^{k+h\phi_{n-1}}$ modulo $p^{n+1}$
go through $g\alpha^k,g\alpha^k+p^n,\ldots,g\alpha^k+(p-1)p^n$. Hence, by
\eqref{equation: Bernoulli relation} in Lemma \ref{lemma: Eg
  relations}, we find
$$
  \pi_n(\div^\infty E_g^{(p^{n+1})})
 =\frac{p^n}2\sum_{k=0}^{\phi_{n-1}-1}B_2\left(
  \frac{g\alpha^k}{p^n}\right)P_{n-1,k}=\div^\infty E_g^{(p^n)}.
$$
When $p|g$, all $ga^{k+h\phi_{n-1}}$ are congruent to $ga^k$ modulo
$p^{n+1}$. Therefore, we have
$$
  \pi_n(\div^\infty E_g^{(p^{n+1})})
 =p^2\cdot\frac{p^n}2\sum_{k=0}^{\phi_{n-1}-1}B_2\left(
  \frac{(g/p)\alpha^k}{p^n}\right)P_{n-1,k}
 =p^2\div^\infty E_{g/p}^{(p^n)}.
$$
This proves the lemma.
\end{proof}

\begin{Lemma} \label{lemma: projection principal to principal} Assume
  $n\ge 1$. If $D$ is a principal divisor in $\DD_n$, then $\pi_n(D)$
  is a principal divisor in $\DD_{n-1}$.

  More precisely, if $f_i$, $i=1,\ldots,\phi_n-1$, is the basis for
  $\FF_n$ given in Theorem \ref{theorem: prime power}, then for
  $p\ge 3$ we have
$$
  \pi_n(\div f_i)=
  \begin{cases}
  0, &i=1,\ldots,\phi_n-\phi_{n-1}, \\
  \div\frac{E^{(p^n)}_{\alpha^{i-1}}(\tau)^{p^2}}
  {E^{(p^n)}_{\alpha^{i+\phi_{n-2}-1}}(\tau)^{p^2}},
  &i=\phi_n-\phi_{n-1}+1,\ldots,\phi_n-\phi_{n-2}, \\
  \qquad \vdots & \qquad \vdots \\
  \div\frac{E^{(p^2)}_{\alpha^{i-1}}(p^{n-2}\tau)^{p^2}}
    {E^{(p^2)}_{\alpha^{i+\phi_0-1}}(p^{n-2}\tau)^{p^2}},
  &i=\phi_n-\phi_1+1,\ldots,\phi_n-\phi_0, \\
  \div\frac{E^{(p)}_{\alpha^{i-1}}(p^{n-1}\tau)^{p^2}}
  {E^{(p)}_{\alpha^i}(p^{n-1}\tau)^{p^2}}, &i=\phi_n-\phi_0+1,\ldots,\phi_n-1.
  \end{cases}
$$
  A similar result also holds for $p=2$.
\end{Lemma}

\begin{proof} Here we prove the case $p$ is an odd prime; the proof of
  the case $p=2$ is similar, and is omitted.

  We first show that $\pi_n(\div f_i)=0$ for
  $i=1,\ldots,\phi_n-\phi_{n-1}$. By Lemma \ref{lemma: divisor of Eg},
  we have
  $$
    \pi_n(\div^\infty E_{\alpha^{i-1}}^{(p^{n+1})})=\div^\infty
    E_{\alpha^{i-1}}^{(p^n)}, \qquad
    \pi_n(\div^\infty E_{\alpha^{i+\phi_{n-1}-1}}^{(p^{n+1})})=\div^\infty
    E_{\alpha^{i+\phi_{n-1}-1}}^{(p^n)}.
  $$
  However, since $\alpha^{\phi_{n-1}}\equiv 1\mod p^n$, we have
  $E_{\alpha^{i-1}}^{(p^n)}=\pm E_{\alpha^{i+\phi_{n-1}-1}}^{(p^n)}$.
  It follows that
  $$
    \pi_n(\div f_i)=\pi_n(\div^\infty f_i)
   =\pi_n(\div^\infty E^{(p^{n+1})}_{\alpha^{i-1}}/
    E_{\alpha^{i+\phi_{n-1}-1}}^{(p^{n+1})})=0
  $$
  for $i=1,\ldots,\phi_n-\phi_{n-1}$.

For $i=\phi_n-\phi_{n-1}+1,\ldots,\phi_n-\phi_{n-2}$, we have, by
\eqref{equation: EN EM},
$$
  f_i=E^{(p^n)}_{\alpha^{i-1}}(p\tau)/
      E^{(p^n)}_{\alpha^{i+\phi_{n-2}-1}}(p\tau)
 =E^{(p^{n+1})}_{p\alpha^{i-1}}(\tau)
 /E^{(p^{n+1})}_{p\alpha^{i+\phi_{n-2}-1}}(\tau).
$$
By Lemma \ref{lemma: divisor of Eg},
$$
  \pi_n(\div f_i)=\pi_n(\div^\infty f_i)=\div^\infty
  \frac{E^{(p^n)}_{\alpha^{i-1}}(\tau)^{p^2}}
  {E^{(p^n)}_{\alpha^{i+\phi_{n-2}-1}}(\tau)^{p^2}}.
$$
Using the criteria given in Theorem \ref{theorem: DQ DO} we find the
last function is in $\FF_{n-1}$ and
$$
  \pi_n(\div f_i)=\div\frac{E^{(p^n)}_{\alpha^{i-1}}(\tau)^{p^2}}
  {E^{(p^n)}_{\alpha^{i+\phi_{n-2}-1}}(\tau)^{p^2}}.
$$
This proves the case $i=\phi_n-\phi_{n-1}+1,\ldots,\phi_n-\phi_{n-2}$.
The remaining cases $i=\phi_n-\phi_{n-2}+1,\ldots,\phi_n-1$ can be
proved in the same way. This gives us the lemma.
\end{proof}

In the next few lemmas, we will establish the fact that
$D\in\DD_{n-1}$ is principal if and only if $\iota_n(D)\in\DD_n$ is
principal.

\begin{Lemma} \label{lemma: iota 1} If $D$ is a principal divisor in
  $\DD_{n-1}$, then $\iota_n(D)$ is a principal divisor in $\DD_n$.
\end{Lemma}

\begin{proof} Let $f^\ast$ be one of the functions in the basis of
  $\FF_{n-1}$ given in Theorem \ref{theorem: prime power} (or Theorem
  \ref{theorem: power of 2} if $p=2$). Define $f(\tau)=f^\ast(p\tau)$.
  From the explicit description of the basis, we see that $f(\tau)$ is
  either one or a product of the functions appearing in our basis for
  $\FF_n$. We now show that $\div f=\iota_n(f^\ast)$.

  Assume $f^\ast(\tau)=\prod_g E_g^{(p^n)}(\tau)^{e_g}$. For a cusp
  $\alpha^k/p^{n+1}\in C_n$, we choose a matrix
  $\sigma=\left(\begin{smallmatrix}\alpha^k&b\\p^{n+1}&d
  \end{smallmatrix}\right)$ in $\Gamma_0(p^{n+1})$. Then we have
  $$
    E_g^{(p^n)}(p\sigma\tau)=E_g^{(p^n)}\left(
    \begin{pmatrix}\alpha^k&pb\\p^n&d\end{pmatrix}(p\tau)\right).
  $$
  Using Lemma \ref{lemma: Eg formulas}, we find
  $$
    E_g^{(p^n)}(p\sigma\tau)=\epsilon E_{\alpha^kg}^{(p^n)}(p\tau)
  $$
  for some root of unity $\epsilon$, and consequently the order of
  $E_g^{(p^n)}(p\tau)$ at $\alpha^k/p^{n+1}$ is
  $$
    p\cdot\frac{p^n}2B_2\left(\frac{\alpha^kg}{p^n}\right),
  $$
  which is the same as $p$ times the order of $E_g^{(p^n)}(\tau)$
  at $\alpha^k/p^n$. From this, we conclude that $\div f=\iota_n(\div f^\ast)$.
  This proves the lemma.
\end{proof}

The proof of the converse statement is more difficult. It relies on
the next two lemmas.

\begin{Lemma} \label{lemma: Bernoulli matrix} Let $N\ge 4$ be an
  integer, $m=\phi(N)/2$, and $a_i$, $1\le i\le m$, be the integers in
  the range $1\le a_i\le N/2$ such that $(a_i,N)=1$. Let $M$ be the
  $m\times m$ matrix whose $(i,j)$-entry is $NB_2(a_ia_j^{-1}/N)/2$,
  where $a_j^{-1}$ denotes the multiplicative inverse of $a_j$ modulo
  $N$. Then we have 
$$
  \det M=\prod_\chi\frac14B_{2,\chi}\neq 0,
$$
where $\chi$ runs over all even characters modulo $N$.
\end{Lemma}

\begin{proof} The proof of
$$
  \det M=\prod_\chi\frac14B_{2,\chi}
$$
can be found in \cite[Lemma 7]{Yang2}, and will not be repeated
here. To see why the determinant is non-zero, we observe that, by a
straightforward computation,
\begin{equation} \label{equation: B0}
  B_{2,\chi_0}=\frac16(1-p)\neq 0,
\end{equation}
where $\chi_0$ is the principal character. Also, Theorem \ref{theorem:
  Yu} in particular implies that
$$
  \prod_{\chi\neq\chi_0\text{ even}}B_{2,\chi}\neq 0.
$$
Therefore, we conclude that $\det M\neq 0$.
\end{proof}

\begin{Lemma} \label{lemma: pre iota 2} Assume $p^{n+1}\ge 5$. Assume
  that $f(\tau)=\prod_g E_g^{(p^{n+1})}(\tau)^{e_g}$ is a modular unit
  in $\FF_n$, where $g\in(\Z/p^{n+1}\Z)^\times/\pm 1$. Suppose that
  for each integer $k$, the orders of $f(\tau)$ at
  $\alpha^{k+h\phi_{n-1}}/p^{n+1}$ take the same values for all
  $h=0,\ldots,p-1$. Then we have $e_g=0$ for all $g$ satisfying
  $p\nmid g$.
\end{Lemma}

\begin{proof} By Lemma \ref{lemma: behavior of Eg}, if $p|g$, then the
  orders of $E_g$ at $\alpha^{k+h\phi_{n-1}}/p^{n+1}$,
  $h=0,\ldots,p-1$, are all $p^{n+1}B_2(\alpha^k(g/p)/p^{n})/2$.
  Therefore, if $f(\tau)=\prod_g E_g^{(p^{n+1})}(\tau)^{e_g}$ has the
  same order at $\alpha^{k+h\phi_{n-1}}/p^{n+1}$ for all
  $h=0,\ldots,p-1$ for any fixed $k$, then the partial product
  $\prod_{p\nmid g}E_g^{e_g}$ also has the same property. Now given
  $k$, let us assume that the order of $\prod_{p\nmid g}E_g^{e_g}$ at
  $\alpha^{k+h\phi_{n-1}}/p^{n+1}$ is $A$. Then we have, by Lemma
  \ref{lemma: behavior of Eg},
  $$
    pA=\sum_{h=0}^{p-1}\sum_{p\nmid g,g\le p^{n+1}/2}e_g\frac{p^{n+1}}2B_2
    \left(\frac{g\alpha^{k+h\phi_{n-1}}}{p^{n+1}}\right).
  $$
  Then by \eqref{equation: Bernoulli relation} in Lemma \ref{lemma: Eg
  relations}, we have
  $$
    pA=\sum_{p\nmid g}e_g\frac{p^n}2B_2\left(\frac{g\alpha^k}{p^n}\right)
    =\sum_{g\le p^n/2,p\nmid g}\frac{p^n}2B_2\left(\frac{g\alpha^k}{p^n}\right)
     \sum_{h=0}^{p-1}e_{g+hp^n}.
  $$
  Now since $f(\tau)$ is assumed to be in $\FF_n$, by Theorem
  \ref{theorem: DQ DO}, we have $\sum_{h=0}^{p-1}e_{g+hp^n}=0$ for all
  $g$. Therefore, we have $A=0$. This is true for all $\alpha^k/p^{n+1}$.
  In other words, we have
  $$
    \sum_{g\le p^{n+1}/2,p\nmid g}e_g
    B_2\left(\frac{g\alpha^k}{p^{n+1}}\right)=0
  $$
  for all $k$. Now write $g=\alpha^j$ and consider the square matrix
  whose $(j,k)$-entry is $B_2(\alpha^{j+k-2}/p^{n+1})$. By Lemma
  \ref{lemma: Bernoulli matrix}, the determinant of this matrix is
  non-zero. Therefore, all $e_g$, $p\nmid g$, are equal to $0$. This
  completes the proof. 
\end{proof}

With the above lemmas, we are now ready to prove the converse to Lemma
\ref{lemma: iota 1}.

\begin{Lemma} \label{lemma: iota 2} Assume that $p$ is a prime and
  $n\ge 1$ is an integer such that $p^n\ge 5$. Let $D$ be
  a divisor in $\DD_{n-1}$. If $\iota_n(D)\in\DD_n$ is principal, then
  $D$ is a principal divisor in $\DD_{n-1}$.
\end{Lemma}

\begin{proof} Let
  $$
    D=\sum_{k=0}^{\phi_{n-1}-1}n_k P_{n-1,k}\in\DD_{n-1}.
  $$
  Assume that $\iota_n(D)$ is principal. That is, assume that there
  exists a function $f(\tau)=\prod_g
  E_g^{(p^{n+1})}(\tau)^{e_g}\in\FF_n$ such that
  $$
    \div f=p\sum_{k=0}^{\phi_{n-1}-1}\sum_{h=0}^{p-1}
    n_kP_{n,k+h\phi_{n-1}}.
  $$
  In other words, we have
  $$
    \frac{p^{n+1}}2\sum_g e_gB_2\left(
    \frac{g\alpha^{k+h\phi_{n-1}}}{p^{n+1}}\right)=pn_k
  $$
  for all $h$ for a given $k$. Since $\iota_n(D)$ has the same order at
  $\alpha^{k+h\phi_{n-1}}/p^{n+1}$ for all $h=0,\ldots,p-1$ for a fixed
  $k$, we have $e_g=0$ whenever $p\nmid g$ by Lemma \ref{lemma: pre
  iota 2}. Thus, we have
  $$
    \frac{p^n}2\sum_{p|g}e_gB_2\left(\frac{(g/p)\alpha^k}{p^n}\right)=n_k,
  $$
  which in turn implies that the function
  $$
    f^\ast(\tau)=\prod_{p|g}E_{g/p}^{(p^n)}(\tau)^{e_g}
  $$
  satisfies $\div f^\ast=D$. It remains to show that $f^\ast$ is a modular
  unit contained in $\FF_{n-1}$, i.e., that $f^\ast$ satisfies
  conditions \eqref{condition DQ} and \eqref{condition DO} of Theorem
  \ref{theorem: DQ DO}.

  Since $f\in\FF_n$, by Theorem \ref{theorem: DQ DO}, the exponents
  $e_g$ satisfy
  $$
    \sum_{g\equiv\pm ap\Mod p^n}e_g=0
  $$
  for all $a$. The same exponents $e_g$ then satisfy
  $$
    \sum_{g:~g/p\equiv\pm a\Mod p^{n-1}}e_g=0,
  $$
  which is condition \eqref{condition DO} for the level $N=p^n$. It
  remains to consider condition \eqref{condition DQ}.

  Observe that $\iota_n(D)$ is a multiple of $p$, whence we have
  \begin{equation} \label{equation: lemma iota 2 1}
    p\,\Big|\sum_{p|g}e_g\frac{p^{n+1}}2
    B_2\left(\frac{g\alpha^k}{p^{n+1}}\right)
  \end{equation}
  for all $k$. We first consider the cases $p\ge 3$. Setting $k=0$ in
  \eqref{equation: lemma iota 2 1}, we have
  $$
    \sum_{p|g}e_g(g^2-gp^{n+1})\equiv 0\mod p^{n+2},
  $$
  or equivalently,
  $$
    \sum_{p|g}e_g(g/p)^2\equiv 0\mod p^n.
  $$
  In other words, $f^\ast$ satisfies the quadratic condition
  \eqref{condition DQ} of Theorem \ref{theorem: DQ DO}. This settles
  the cases $p\ge 3$.

  For $p=2$, \eqref{equation: lemma iota 2 1} with $k=0$ yields
  $$
    \sum_{2|g}e_g(g^2-2^{n+1}g)\equiv 0\mod 2^{n+3},
  $$
  i.e.,
  $$
    \sum_{2|g}e_g\left((g/2)^2-2^n(g/2)\right)\equiv 0\mod 2^{n+1}.
  $$
  Partition the sum $\sum_ge_g(g/2)$ into two parts $\sum_{g\equiv 0\Mod
  4}$ and $\sum_{g\equiv 2\Mod 4}$. For the terms with $4|g$, we
  clearly have
  $$
    \sum_{g\equiv 0\Mod 4}e_g(g/2)\equiv 0\mod 2.
  $$
  For the terms with $g\equiv 2\Mod 4$, we have
  $$
    \sum_{g\equiv 2\Mod 4}e_g(g/2)\equiv\sum_{g\equiv 2\Mod 4}e_g\mod 2.
  $$
  Since $e_g$ satisfy condition \eqref{condition DO} for $N=2^{n+1}$, we
  must have
  $$
    \sum_{g\equiv 2\Mod 4}e_g=0.
  $$
  Therefore,
  $$
    \sum_{2|g}e_g(g/2)\equiv 0\mod 2.
  $$
  It follows that
  $$
    \sum_{2|g}e_g(g/2)^2\equiv
    \sum_{2|g}e_g\left((g/2)^2-2^n(g/2)\right)\equiv 0\mod 2^{n+1},
  $$
  which is \eqref{condition DQ} for $N=2^n$. This proves the case
  $p=2$, and the proof of the lemma is complete.
\end{proof}

From Lemmas \ref{lemma: iota 1} and \ref{lemma: iota 2}, we
immediately get the following corollary.

\begin{Corollary} The homomorphism $\iota_n$ induces an embedding
  $\iota_n^\ast:\CC_{n-1}\to\CC_n$ given by
  $\iota_n^\ast([D])=[\iota_n(D)]$.
\end{Corollary}
\end{section}

%%%%%%%%%%%%%%%%%%%%%%%%%%%%%%%%%%%%%%%%%%%%%%%%%%%%%%%%%%%%
%
%  Propositions
%
%%%%%%%%%%%%%%%%%%%%%%%%%%%%%%%%%%%%%%%%%%%%%%%%%%%%%%%%%%%%

\begin{section}{Proof of Propositions}
\label{section: propositions}
\begin{subsection}{Proof of Proposition \ref{proposition: C0}}
\label{subsection: C0}
\begin{Lemma} \label{lemma: Bernoulli congruence} Let $p\ge 3$ be an
  odd prime. Let $\omega$ be a generator of the group of Dirichlet
  characters modulo $p$. Then we have the congruence
  $$
    p\prod_{i=1}^{(p-1)/2-1}B_{2,\omega^{2i}}\equiv
    \begin{cases}\displaystyle
    -\prod_{i=1}^{(p-1)/2-2}\frac{B_{2i+2}}{i+1}\mod p,
    &\text{if }p\ge 5, \\
    -1 \mod 3, &\text{if }p=3, \end{cases}
  $$
  where $B_{2,\omega^{2i}}$ are the generalized Bernoulli numbers and
  $B_{2i+2}$ are Bernoulli numbers.
\end{Lemma}

\begin{proof} The case $p=3$ can be verified directly. We now assume
  that $p\ge 5$.

  Since the product is a rational number, we may regard $\omega$ as
  the Teichm\"uller character $\omega:~\Z_p^\times\to\mu_{p-1}$ from
  $\Z_p^\times$ to the group of $(p-1)$-st roots of unity in $\Z_p$
  characterized by $\omega(a)\equiv a\mod p$ for all
  $a\in\Z_p^\times$. It is well-known that for $2i\neq p-3$,
  $B_{2,\omega^{2i}}$ is contained in $\Z_p$ and satisfies
  $$
    B_{2,\omega^{2i}}\equiv\frac{B_{2i+2}}{i+1} \mod p.
  $$
  (For a proof, follow the argument in \cite[Corollary
  5.15]{Washington}.) Also, for $2i=p-3$, we have
  $$
    pB_{2,\omega^{p-3}}=\sum_{a=1}^{p-1}\omega^{-2}(a)(a^2-pa+p^2/6)
    \equiv\sum_{a=1}^{p-1}\omega^{-2}(a)a^2\equiv\sum_{a=1}^{p-1}1
    \equiv-1\mod p.
  $$
  Then the lemma follows.
\end{proof}

\begin{proof}[Proof of Proposition \ref{proposition: C0}]
  The cases $p=2,3$ can be easily seen from the fact that the modular
  curves $X_1(8)$ and $X_1(9)$ have genus zero. Now assume $p\ge 5$.
  By Theorem \ref{theorem: Yu}, the order of the divisor group $\CC_0$
  is
  $$
    |\CC_0|=p\prod_{i=1}^{(p-3)/2}\frac14B_{2,\omega^{2i}}.
  $$
  Using Lemma \ref{lemma: Bernoulli congruence}, we obtain
  $$
    |\CC_0|\equiv-\frac14\prod_{i=1}^{(p-5)/2}
    \frac14\frac{B_{2i+2}}{i+1}\mod p.
  $$
  By the assumption that $p$ is a regular prime, none of
  $B_4,\ldots,B_{p-3}$ is divisible by $p$. Therefore, $p$ does not
  divide $|\CC_0|$.
\end{proof}
\end{subsection}

\begin{subsection}{Proof of Proposition \ref{proposition: p-rank}}
\label{subsection: p-rank}
  Among the five propositions, this proposition is perhaps the most
  complicated to prove.

  Recall that given a free $\Z$-module $\Lambda$ of finite rank $r$
  with basis $\{a_1,\ldots,a_r\}$ and a submodule $\Lambda'$ generated
  by $b_1,\ldots,b_s$ with $b_i=\sum r_{ij}a_j$, the standard method
  to determine the group structure of $\Lambda/\Lambda'$ is to compute
  the Smith normal form of the matrix $(r_{ij})$. Then the $p$-rank of
  the group $\Lambda/\Lambda'$ is simply the number of diagonals in
  the Smith normal form that are divisible by $p$. Thus, in order to
  prove Proposition \ref{proposition: p-rank}, we need to know very
  precisely the linear dependence over $\F_p$ among the divisors of
  modular units generating $\FF_n$. In the first two lemmas, we will
  show that the divisors of the first $\phi_n-\phi_{n-1}$ functions in
  the basis for $\FF_n$ are linearly independent over $\F_p$.

\begin{Lemma} \label{lemma: determinant of new part} Let $p$ be a
  prime and $n\ge 1$ be an integer such that $p^{n+1}\ge 5$. Let
  $\alpha$ be a generator of $(\Z/p^{n+1}\Z)^\times/\pm 1$. Let $f_i$,
  $i=1,\ldots,\phi_n-1$, be the basis for $\FF_n$ given in
  Theorem \ref{theorem: prime power} or Theorem \ref{theorem: power of
  2}. Let $M=(m_{ij})$ be the square matrix of size
  $\phi_n-\phi_{n-1}$ such that $m_{ij}$ is the order of
  $f_i$ at the cusp $\alpha^{j-1}/p^{n+1}$. Then we have
  $$
    \det M=\epsilon p\prod_{\chi\text{ even primitive}}\frac14B_{2,\chi},
  $$
  where $\chi$ runs over all even primitive Dirichlet characters
  modulo $p^{n+1}$ and $\epsilon$ is either $1$ or $-1$.
\end{Lemma}

\begin{proof} Let $A=(a_{ij})$ be the $\phi_n\times\phi_n$ matrix with
  $a_{ij}=p^{n+1}B_2(\alpha^{i+j-2}/p^{n+1})/2$, which is the order of
  $E_{\alpha^{i-1}}$ at $P_{n,j-1}=\alpha^{j-1}/p^{n+1}$. Define
  $$
    V_1=\begin{pmatrix}
    I & -I & 0 & \cdots & \cdots & \cdots \\
    0 & I & -I & \cdots & \cdots & \cdots \\
    \vdots & \vdots & & & \vdots & \vdots \\
    \cdots & \cdots & \cdots & I & -I & 0 \\
    \cdots & \cdots & \cdots & 0 & I & -I \\
    I & I & \cdots & \cdots & I & I \end{pmatrix},
  $$
  where the matrix consists of $p^2$ blocks, each of which is of size
  $\phi_{n-1}\times\phi_{n-1}$, and $I$ is the identity matrix of
  dimension $\phi_{n-1}$.
  Let $\beta$ be an integer such that $\alpha\beta\equiv 1\mod p$.
  Set also
  $$
    V_2=\begin{pmatrix}
    1 & -\beta^2 & 0 & \cdots & \cdots & \cdots \\
    0 & 1 & -\beta^2 & \cdots & \cdots & \cdots \\
    \vdots & \vdots & & & \vdots & \vdots \\
    \cdots & \cdots & \cdots & 1 & -\beta^2 & 0 \\
    \cdots & \cdots & \cdots & 0 & p & 0 \\
    0 & 0 & \cdots & \cdots & 0 & I\end{pmatrix},
  $$
  where the identity matrix at the lower right corner has dimension
  $\phi_{n-1}$. Then for $i=1,\ldots,\phi_n-\phi_{n-1}$, the
  $(i,j)$-entry of the matrix $V_2V_1A$ is the order of $f_i$ at
  $P_{n,j-1}$, while for $i=\phi_n-\phi_{n-1}+1,\ldots,\phi_n$, the
  $(i,j)$-entry of $V_2V_1A$ is
  $$
    \frac{p^{n+1}}2\sum_{h=0}^{p-1}
    B_2\left(\frac{\alpha^{i+j+h\phi_{n-1}-2}}{p^{n+1}}\right).
  $$
  By \eqref{equation: Bernoulli relation} in Lemma
  \ref{lemma: Eg relations}, this is equal to
  \begin{equation} \label{equation: A'}
    \frac{p^n}2B_2\left(\frac{\alpha^{i+j-2}}{p^n}\right).
  \end{equation}
  Observe that
  $B_2(\alpha^{i+j-2}/p^n)=B_2(\alpha^{i+j+k\phi_{n-1}-2}/p^n)$ for
  all integers $k$. That is, $V_2V_1A$ takes the form
  $$
    V_2V_1A=\begin{pmatrix}
      \\
    \text{ order of }f_i\text{ at }\alpha^{j-1}/p^{n+1} \\
    \text{ for }i=1,\ldots,\phi_n-\phi_{n-1} \\
      \\
    A' \quad A' \quad \cdots \quad \cdots \quad A'\quad A'
    \end{pmatrix},
  $$
  where $A'$ is a square matrix of size $\phi_{n-1}$ whose
  $(i,j)$-entry is given by \eqref{equation: A'}.

  Now let
  $$
    U_1=\begin{pmatrix}
    I & 0 & \cdots & \cdots & 0 & I \\
    0 & I & \cdots & \cdots & 0 & I \\
    \vdots & \vdots & & & \vdots & \vdots \\
    \vdots & \vdots & & & \vdots & \vdots \\
    0 & 0 & \cdots & \cdots & I & I \\
    0 & 0 & \cdots & \cdots & 0 & I \end{pmatrix},
  $$
  and consider $V_2V_1AU_1$. For $i=1,\ldots,\phi_n-\phi_{n-1}$ and
  $j=\phi_n-\phi_{n-1}+1,\ldots,\phi_n$, the $(i,j)$-entry of
  $V_2V_1AU_1$ is
  $$
    \sum_{h=0}^{p-1}(\text{order of }f_i\text{ at }
    P_{n,j+h\phi_{n-1}-1}).
  $$
  By Lemma \ref{lemma: projection principal to principal}, this sum is
  equal to $0$. In other words,
  $$
    V_2V_1AU_1=\begin{pmatrix}
     & & & & & 0 \\
     & & & & & 0 \\
     & & M & & & \vdots \\
     & & & & & 0 \\
     & & & & & 0 \\
    A' & \cdots & \cdots & \cdots & A' & pA'\end{pmatrix},
  $$
  where $M$ is the $(\phi_n-\phi_{n-1})\times(\phi_n-\phi_{n-1})$
  matrix specified in the lemma. This shows that
  $$
    \det(V_2V_1AU_1)=p^{\phi_{n-1}}(\det M)(\det A').
  $$
  On the other hand, we have, by Lemma \ref{lemma: Bernoulli matrix},
  $$
    \det A=\pm\prod_{\chi\Mod p^{n+1}}\frac14B_{2,\chi}, \qquad
    \det A'=\pm\prod_{\chi\Mod p^n}\frac14B_{2,\chi}.
  $$
  Also,
  $$
    \det V_1=p^{\phi_{n-1}}, \qquad \det V_2=p, \qquad
    \det U_1=1.
  $$
  Combining everything, we conclude that
  $$
    \det M=\pm p\prod_{\chi\Mod p^{n+1}}\frac14B_{2,\chi}\Big/
    \prod_{\chi\Mod p^n}\frac14B_{2,\chi}
  =\pm p\prod_{\chi\text{ even primitive mod }p^{n+1}}\frac14
    B_{2,\chi},
  $$
  as claimed in the lemma.
\end{proof}

Here we give an example to exemplify the above argument.

\begin{Example} \upshape Consider the case $p=3$ and $n=2$. We choose
  $\alpha=2$ and $\beta=-1$. With the notations as above, we have
$$
  A=\frac1{108}\left(\begin{smallmatrix}
  191&143&59&-61&-109&23&-97&-37&-121 \\
  143&59&-61&-109&23&-97&-37&-121&191 \\
  59&-61&-109&23&-97&-37&-121&191&143 \\
  -61&-109&23&-97&-37&-121&191&143&59 \\
  -109&23&-97&-37&-121&191&143&59&-61 \\
  23&-97&-37&-121&191&143&59&-61&-109 \\
  -97&-37&-121&191&143&59&-61&-109&23 \\
  -37&-121&191&143&59&-61&-109&23&-97 \\
  -121&191&143&59&-61&-109&23&-97&-37
  \end{smallmatrix}\right),
$$
where the $(i,j)$-entry is $27B_2(2^{i+j-2}/27)/2$,
$$
  V_2=\left(\begin{smallmatrix}
  1&-1&0&0&0&0&0&0&0 \\
  0&1&-1&0&0&0&0&0&0 \\
  0&0&1&-1&0&0&0&0&0 \\
  0&0&0&1&-1&0&0&0&0 \\
  0&0&0&0&1&-1&0&0&0 \\
  0&0&0&0&0&3&0&0&0 \\
  0&0&0&0&0&0&1&0&0 \\
  0&0&0&0&0&0&0&1&0 \\
  0&0&0&0&0&0&0&0&1\end{smallmatrix}\right), \qquad
  V_1=\left(\begin{smallmatrix}
  1&0&0&-1&0&0&0&0&0 \\
  0&1&0&0&-1&0&0&0&0 \\
  0&0&1&0&0&-1&0&0&0 \\
  0&0&0&1&0&0&-1&0&0 \\
  0&0&0&0&1&0&0&-1&0 \\
  0&0&0&0&0&1&0&0&-1 \\
  1&0&0&1&0&0&1&0&0 \\
  0&1&0&0&1&0&0&1&0 \\
  0&0&1&0&0&1&0&0&1\end{smallmatrix}\right).
$$
Then
$$
  V_2V_1A=\left(\begin{smallmatrix}
  0&2&0&1&-2&4&-1&0&-4 \\
  2&0&1&-2&4&-1&0&-4&0 \\
  0&1&-2&4&-1&0&-4&0&2 \\
  1&-2&4&-1&0&-4&0&2&0 \\
  -2&4&-1&0&-4&0&2&0&1 \\
  4&-8&-5&-5&7&7&1&1&-2\\
  a&b&c&a&b&c&a&b&c \\
  b&c&a&b&c&a&b&c&a \\
  c&a&b&c&a&b&c&a&b \end{smallmatrix}\right), \quad
  \begin{cases}
  a=11/36=9B_2(1/9)/2, \\
  b=-1/36=9B_2(2/9)/2, \\
  c=-13/36=9B_2(4/9)/2. \end{cases}
$$
Here the first $6$ rows are the orders of
$$
  \frac{E_1E_{11}}{E_2E_8},\ \frac{E_2E_5}{E_4E_{11}},\
  \frac{E_4E_{10}}{E_8E_5},\ \frac{E_8E_7}{E_{11}E_{10}},\
  \frac{E_{11}E_{13}}{E_5E_7},\ \frac{E_5^3}{E_{13}^3}
$$
at the cusps $2^{j-1}/27$. The matrices $U_1$ and $V_2V_1AU_1$ then are
$$
  U_1=\left(\begin{smallmatrix}
  1&0&0&0&0&0&1&0&0 \\
  0&1&0&0&0&0&0&1&0 \\
  0&0&1&0&0&0&0&0&1 \\
  0&0&0&1&0&0&1&0&0 \\
  0&0&0&0&1&0&0&1&0 \\
  0&0&0&0&0&1&0&0&1 \\
  0&0&0&0&0&0&1&0&0 \\
  0&0&0&0&0&0&0&1&0 \\
  0&0&0&0&0&0&0&0&1\end{smallmatrix}\right), \quad
  V_2V_1AU_1=\left(\begin{smallmatrix}
  0&2&0&1&-2&4&0&0&0 \\
  2&0&1&-2&4&-1&0&0&0 \\
  0&1&-2&4&-1&0&0&0&0 \\
  1&-2&4&-1&0&-4&0&0&0 \\
  -2&4&-1&0&-4&0&0&0&0 \\
  4&-8&-5&-5&7&7&0&0&0\\
  a&b&c&a&b&c&3a&3b&3c \\
  b&c&a&b&c&a&3b&3c&3a \\
  c&a&b&c&a&b&3c&3a&3b \end{smallmatrix}\right).
$$
We find
$$
  \det M=\det\left(\begin{smallmatrix}
  0&2&0&1&-2&4 \\
  2&0&1&-2&4&-1 \\
  0&1&-2&4&-1&0& \\
  1&-2&4&-1&0&-4 \\
  -2&4&-1&0&-4&0 \\
  4&-8&-5&-5&7&7 \end{smallmatrix}\right)=-5833
 =-3\prod_{\chi\text{ even primitive mod }27}\frac14B_{2,\chi}.
$$
\end{Example}

\begin{Lemma} \label{lemma: Bernoulli product} Let $p$ be a regular
  prime and $n\ge 1$ be an integer. Then we have 
  $$
    p\prod_{\chi\text{ even primitive}}\frac14B_{2,\chi}\equiv 1\mod p,
  $$
  where the product runs over all even primitive Dirichlet characters
  modulo $p^{n+1}$.
\end{Lemma}

\begin{proof} First of all, for any non-trivial even Dirichlet character
  $\chi$ we have
  $$
    \sum_{a=1}^{p^{n+1}}\chi(a)=0
  $$
  and
  $$
    \sum_{a=1}^{p^{n+1}}a\chi(a)=\frac12\sum_{a=1}^{p^{n+1}}
    \big(a\chi(a)+(p^{n+1}-a)\chi(p^{n+1}-a)\big)
   =\frac{p^{n+1}}2\sum_{a=1}^{p^{n+1}}\chi(a)=0.
  $$
  Thus,
  \begin{equation} \label{equation: lemma Bernoulli product 1}
    B_{2,\chi}=p^{n+1}\sum_{a=1}^{p^{n+1}}\left(
    \frac{a^2}{p^{2n+2}}-\frac a{p^{n+1}}+\frac16\right)\chi(a)
   =\frac1{p^{n+1}}\sum_{a=1}^{p^{n+1}}\chi(a)a^2.
  \end{equation}
  Now we consider the case $p$ is an odd regular prime first.

  Fix a generator $\alpha$ of the multiplicative group
  $(\Z/p^{n+1}\Z)^\times$. For a non-negative integer $m$, write
  $r(m)=\lfloor\alpha^m/p^{n+1}\rfloor$ and
  $s(m)=\alpha^m/p^{n+1}-r(m)$. We have
  $$
    p^{n+1}s(m)^2=\frac{\alpha^{2m}}{p^{n+1}}-2\alpha^m r(m)+
    p^{n+1}r(m)^2
  $$
  Therefore, if $a$ is the integer in the range $0<a<p^{n+1}$ such that
  $\alpha^m\equiv a\Mod p^{n+1}$, then
  \begin{equation} \label{equation: lemma Bernoulli product 4}
    \frac{a^2}{p^{n+1}}-\frac{\alpha^{2m}}{p^{n+1}}
   =-2\alpha^m r(m)+p^{n+1}r(m)^2
  \end{equation}
  is an integer. Denote this integer by $\delta(m)$. Then by
  \eqref{equation: lemma Bernoulli product 1}, we may write
  \begin{eqnarray} \nonumber
    B_{2,\chi}&=&\frac1{p^{n+1}}\sum_{a=1}^{p^{n+1}}\chi(a)a^2 \\
   &=&\frac1{p^{n+1}}\sum_{m=0}^{p^n(p-1)-1}\chi(\alpha^m)\alpha^{2m}
   +\sum_{m=0}^{p^n(p-1)-1}\chi(\alpha^m)\delta(m) \nonumber \\
  &=&\frac{1-\alpha^{2p^n(p-1)}}{p^{n+1}(1-\chi(\alpha)\alpha^2)}
  +\sum_{m=0}^{p^n(p-1)-1}
    \chi(\alpha^m)\delta(m). \label{equation: lemma Bernoulli product 2}
  \end{eqnarray}
  Note that the number $(1-\alpha^{2p^n(p-1)})/p^{n+1}$ is an integer.
  Therefore, $(1-\chi(\alpha)\alpha^2)B_{2,\chi}$ is an algebraic integer.

  Let $\omega$ and $\theta$ denote the Dirichlet characters satisfying
  $$
    \omega(\alpha)=\zeta_{p-1}, \qquad
    \theta(\alpha)=\zeta_{p^n},
  $$
  respectively, where $\zeta_m=e^{2\pi i/m}$.
  Set $\chi_{ij} = \omega^{2i}\theta^j$.
  Then the set of even
  primitive Dirichlet character modulo $p^{n+1}$ is precisely
  $$
    \{\chi_{ij}=\omega^{2i}\theta^j:~0\le i<(p-1)/2,~0\le j<p^n,
    ~p\nmid j\}.
  $$
  From \eqref{equation: lemma Bernoulli product 2}, we have, for all
  $j$ not divisible by $p$,
  \begin{equation*}
  \begin{split}
  &(1-\omega^{2i}(\alpha)\alpha^2)B_{2,\omega^{2i}}
  -(1-\chi_{ij}(\alpha)\alpha^2)B_{2,\chi_{ij}} \\
 =\ &(1-\omega^{2i}(\alpha)\alpha^2)\sum_{m=0}^{p^n(p-1)-1}
   \omega^{2i}(\alpha^m)\delta(m)
   -(1-\chi_{ij}(\alpha)\alpha^2)
   \sum_{m=0}^{p^n(p-1)-1}\chi_{ij}(\alpha^m)\delta(m) \\
 =\ &(1-\omega^{2i}(\alpha)\alpha^2)\sum_{m=0}^{p^n(p-1)-1}
  \omega^{2i}(\alpha^m)\delta(m)(1-\theta^j(\alpha^m)) \\
  &\qquad\quad-\omega^{2i}(\alpha)\alpha^2(1-\theta^j(\alpha))
   \sum_{m=0}^{p^n(p-1)-1}\chi_{ij}(\alpha^m)\delta(m) \\
 \equiv\ &0\mod 1-\zeta_{p^n}.
  \end{split}
  \end{equation*}
  (Note that when $i=0$, $\omega^0=\chi_0$ is principal, and
  \eqref{equation: lemma Bernoulli product 1} does not hold in this
  case. However, the difference is $p^n$ times a $p$-unit, and the
  above congruence still holds.) In other words,
  $$
    \prod_{j=1,\,p\nmid j}^{p^n}(1-\chi_{ij}(\alpha)\alpha^2)
      B_{2,\chi_{ij}}\equiv
    \big((1-\omega^{2i}(\alpha)\alpha^2)B_{2,\omega^{2i}}\big)^{p^{n-1}(p-1)}
    \mod 1-\zeta_{p^n}.
  $$
  It follows that
  \begin{equation*}
  \begin{split}
  &\prod_{\chi\text{ even primitive}}(1-\chi(\alpha)\alpha^2)B_{2,\chi}
   =\prod_{i=0}^{(p-1)/2-1}\prod_{j=1,\,p\nmid j}^{p^n}
    (1-\chi_{ij}(\alpha)\alpha^2)B_{2,\chi_{ij}} \\
  &\qquad\qquad\equiv\left(\prod_{i=0}^{(p-1)/2-1}
   (1-\omega^{2i}(\alpha)\alpha^2)
    B_{2,\omega^{2i}}\right)^{p^{n-1}(p-1)}\mod 1-\zeta_{p^n}.
  \end{split}
  \end{equation*}
  Now consider the product in the last expression. We have
  $$
    \prod_{i=0}^{(p-1)/2-1}(1-\omega^{2i}(\alpha)\alpha^2)=1-\alpha^{p-1}.
  $$
  Since $\alpha$ is a generator for $(\Z/p^n\Z)^\times$ for all $n$, we
  have $1-\alpha^{p-1}=pu$ for some integer $u$ relatively prime to $p$.
  Also, according to \eqref{equation: B0} and Lemma
  \ref{lemma: Bernoulli congruence}, we have
  $$
   p\prod_{i=0}^{(p-1)/2-1}B_{2,\omega^{2i}}\equiv
   \begin{cases}\displaystyle-\frac16\prod_{i=1}^{(p-1)/2-2}
     \frac{B_{2i+2}}{i+1}\mod p, &\text{if }p\ge 5, \\
   -1, &\text{if }p=3. \end{cases}
  $$
  By the assumption that $p$ is a regular prime, this product is
  relatively prime to $p$. Therefore, we have
  $$
    \left(\prod_{i=0}^{(p-1)/2-1}(1-\omega^{2i}(\alpha)\alpha^2)
    B_{2,\omega^{2i}}\right)^{p-1}\equiv 1\mod p,
  $$
  and consequently
  $$
    \prod_{\chi\text{ even primitive}}(1-\chi(\alpha)\alpha^2)B_{2,\chi}
    \equiv 1\mod 1-\zeta_{p^n}.
  $$
  Since the product is a rational integer, the congruence actually
  holds modulo $p$. Finally, because $\alpha$ is a generator of
  $(\Z/p^n\Z)^\times$ for all $n$, there exists an integer $u$
  relatively prime to $p$ such that $\alpha^{p^k(p-1)}\equiv
  1-up^{k+1}\mod p^{k+2}$ for all $k\ge 0$. 
  Thus,
  \begin{equation} \label{lemma: Bernoulli congruence 5}
    \prod_{\chi\text{ even primitive}}(1-\chi(\alpha)\alpha^2)
   =\frac{1-\alpha^{p^n(p-1)}}{1-\alpha^{p^{n-1}(p-1)}}
   =\frac{up^{n+1}+\cdots}{up^n+\cdots}
   \equiv p\mod p^2.
  \end{equation}
  From this we conclude that
  $$
    p\prod_{\chi\text{ even primitive}}\frac14B_{2,\chi}\equiv 1\mod p.
  $$
  This completes the proof of the case $p$ is an odd regular prime.

  Now consider the case $p=2$ with $n\ge 2$. Choose $\alpha=3$ to be a
  generator of $(\Z/2^{n+1}\Z)^\times/\pm 1$. Set $\zeta=e^{2\pi
  i/2^{n-1}}$, and let $\theta$ be the Dirichlet character satisfying
  $\theta(-1)=1$ and $\theta(3)=\zeta$. Then the set of even primitive
  Dirichlet characters modulo $2^{n+1}$ is
  $$
    \{\theta^j:~1\le j\le 2^{n-1},~2\nmid j\}.
  $$
  Since $\theta$ is even, we have
  $$
    \frac14B_{2,\theta^j}=\frac{2^{n+1}}2
    \sum_{a\in(\Z/2^{n+1}\Z)^\times/\pm 1}\theta^j(a)
    B_2\left(\frac a{2^{n+1}}\right).
  $$
  By a similar calculation as before, we find that if $\theta^j$ is not
  principal, then
  $$
    \frac14B_{2,\theta^j}=\frac1{2^{n+2}}
    \sum_{a\in(\Z/2^{n+1}\Z)^\times/\pm 1}\theta^j(a)a^2.
  $$
  Now for a non-negative integer $m$, define
  $$
    \delta(m)=-\frac{3^{2m}}{2^{n+1}}+2^{n+1}\left\{\frac{3^m}{2^{n+1}}
    \right\}^2
  $$
  as in \eqref{equation: lemma Bernoulli product 4}. Following the
  computation in \eqref{equation: lemma Bernoulli product 2}, we get
  $$
    \frac14B_{2,\theta^j}=\frac{1-3^{2^n}}{2^{n+2}(1-9\zeta^j)}
   +\frac12\sum_{m=0}^{\phi(2^{n+1})/2-1}\theta^j(3^m)\delta(m).
  $$
  Now we have $3^{2^n}=(1+8)^{2^{n-1}}\equiv 1+2^{n+2}\Mod 2^{n+3}$.
  Also, from \eqref{equation: lemma Bernoulli product 4}, we see that
  $\delta(m)$ is always even. Thus, $(1-9\zeta^j)B_{2,\theta^j}/4$ is
  an algebraic integer. By the same argument as before, we find
  $$
    \frac{1-9}4B_{2,\chi_0}-\frac{1-9\zeta^j}4B_{2,\theta^j}
    \equiv 0\mod 1-\zeta,
  $$
  for all odd $j$ and thus
  $$
    \prod_{\chi\text{ even primitive}}\frac{1-9\chi(3)}4B_{2,\chi}
    \equiv 1\mod 2.
  $$
  Finally, as \eqref{lemma: Bernoulli congruence 5}, we have
  $$
    \prod_{\chi\text{ even primitive}}(1-9\chi(3))\equiv 2\mod 4.
  $$
  This proves the case $p=2$.
\end{proof}

\begin{proof}[Proof of Proposition \ref{proposition: p-rank}]
  Let $\alpha$ be a generator of $(\Z/p^{n+1}\Z)^\times/\pm 1$.
  Specifically, for $p=2$, we set $\alpha=3$, and for an odd prime
  $p$, we let $\alpha$ be an integer such that $\alpha$ generates
  $(\Z/p\Z)^\times$, but $\alpha^{p-1}\not\equiv 1\mod p^2$. Let
  $f_i$, $i=1,\ldots,\phi_n-1$, be the generators of $\FF_n$ given in
  Theorem \ref{theorem: prime power} or Theorem
  \ref{theorem: power of 2}. Let $M$ be the $(\phi_n-1)\times\phi_n$
  matrix whose $(i,j)$-entry is the order of $f_i$ at
  $\alpha^{j-1}/p^{n+1}$. Let $U$ and $V$ be the unimodular matrices
  such that $M'=UMV$ is in the Smith normal form. That is, if
  $M'=(m_{ij})$, then
  \begin{enumerate}
  \item $m_{11}|m_{22}|\cdots$, and
  \item $m_{ij}=0$ if $i\neq j$.
  \end{enumerate}
  ($m_{ii}\neq 0$ for all $i$ since the rank of $M$ is $\phi_n-1$.)
  Then the $p$-rank of $\CC_n$ is equal to the number of $m_{ii}$ that
  are divisible by $p$. In other words, if we consider $M$ as a matrix
  over $\F_p$, then our $p$-rank is actually equal to
  $$
    \phi_n-1-(\text{the rank of }M\text{ over }\F_p).
  $$
  We now determine the rank of $M$ over $\F_p$.

  From Lemmas \ref{lemma: determinant of new part} and \ref{lemma:
  Bernoulli product}, we know that the first $\phi_n-\phi_{n-1}$ rows
  of $M$ are linearly independent over $\F_p$. Thus, the rank of $M$
  over $\F_p$ is at least $\phi_n-\phi_{n-1}=p^{n-1}(p-1)^2/2$. It
  remains to prove that the remaining rows are all linearly dependent
  of the first $\phi_n-\phi_{n-1}$ rows modulo $p$.

  We first consider row $\phi_n-\phi_{n-1}+1$ to row $\phi_n-\phi_0$.
  (For $p=2$, consider row $\phi_n-\phi_{n-1}+1$ to row $\phi_n-\phi_2$.)
  Let $\ell$ be an integer between $1$ and $n-1$. (For $p=2$, let
  $1\le\ell\le n-3$.) By Theorems
  \ref{theorem: prime power} and \ref{theorem: power of 2}, for $i$
  from $\phi_n-\phi_{n-\ell}+1$ to $\phi_n-\phi_{n-\ell-1}$, the $i$th
  row of $M$ is the divisor of the function
  $$
    f_i=E^{(p^{n-\ell+1})}_{\alpha^{i-1}}(p^\ell\tau)/
        E^{(p^{n-\ell+1})}_{\alpha^{i+\phi_{n-\ell-1}-1}}(p^\ell\tau),
  $$
  which by Lemma \ref{lemma: behavior of Eg}, is
  \begin{equation} \label{equation: proposition 2 1}
    p^\ell\cdot\frac{p^{n-\ell+1}}2\sum_{k=0}^{\phi_n-1}\left(
    B_2\left(\frac{\alpha^{i+k-1}}{p^{n-\ell+1}}\right)-
    B_2\left(\frac{\alpha^{i+\phi_{n-\ell-1}+k-1}}{p^{n-\ell+1}}\right)
    \right)P_{n,k}.
  \end{equation}
  Now $\alpha^{\phi_{n-\ell-1}}\equiv-(1+up^{n-\ell})\mod
  p^{n-\ell+1}$ for some integer $u$ not divisible by $p$. (For $p=2$,
  we have $\alpha^{\phi_{n-\ell-1}}\equiv 1+2^{n-\ell}\mod
  2^{n-\ell+1}$ instead when $n-\ell\ge 3$.) Then a straightforward
  calculation gives
  \begin{equation*}
  \begin{split}
   &\frac{p^{n-\ell+1}}2\left(B_2\left(
    \frac{\alpha^{i+k-1}}{p^{n-\ell+1}}\right)
   -B_2\left(\frac{\alpha^{i+\phi_{n-\ell-1}+k-1}}{p^{n-\ell+1}}\right)
    \right)
%    \equiv\frac{\alpha^{i+k-1}(1-\alpha^{\phi_{n-\ell-1}})}{p^{n-\ell+1}}
    \equiv-\frac{u\alpha^{2(i+k-1)}}p\mod 1.
  \end{split}
  \end{equation*}
  This shows that if $\ell\ge 2$, then the divisor of $f_i$ for
  $i$ from $\phi_n-\phi_{n-\ell}+1$ to $\phi_n-\phi_{n-\ell-1}$ is
  divisible by $p$. For such $\ell$, the rows do not contribute
  anything to the rank of $M$ over $\F_p$.

  When $\ell=1$, the above computation shows that the $i$th row of $M$
  for $i$ from $\phi_n-\phi_{n-1}+1$ to $\phi_n-\phi_{n-2}$ is
  congruent to
  $$
    -u\alpha^{2(i-1)}(1,\alpha^2,\alpha^4,\ldots,\alpha^{2\phi_n-2})
  $$
  modulo $p$. On the other hand, the $(\phi_n-\phi_{n-1})$-th row of
  $M$ is the divisor of
  $$
    E^p_{\alpha^{\phi_n-\phi_{n-1}-1}}/E^p_{\alpha^{\phi_n-1}}.
  $$
  By a similar computation, we find that it is congruent to
  $$
    -u\alpha^{2(\phi_n-\phi_{n-1}-1)}(1,\alpha^2,\alpha^4,
     \ldots,\alpha^{2\phi_n-2}).
  $$
  From this we see that row $\phi_n-\phi_{n-1}+1$ to row
  $\phi_n-\phi_0$ of $M$ are all multiples of the
  $(\phi_n-\phi_{n-1})$-th row of $M$ modulo $p$.

  Finally, for $i=\phi_n-\phi_0+1,\ldots,\phi_n-1$, we find that the
  $i$th row is congruent to
  $$
    (1-\alpha^2)\frac{\alpha^{2i-2}}2
    (1,\alpha^2,\alpha^4,\ldots,\alpha^{2\phi_n-2})
  $$
  modulo $p$, which again is a multiple of the
  $(\phi_n-\phi_{n-1})$-th row of $M$ modulo $p$. Therefore, the rank
  of $M$ over $\F_p$ is precisely $\phi_n-\phi_{n-1}$. We conclude
  that the $p$-rank of $\CC_n$ is
  $$
    \phi_n-1-(\phi_n-\phi_{n-1})=\phi_{n-1}-1=p^{n-1}(p-1)/2-1.
  $$
  This completes the proof of the proposition.
\end{proof}
\end{subsection}

\begin{subsection}{Proof of Proposition \ref{proposition: index of
  projection}}
\label{subsection: index of projection}
Let $f_i$,
  $i=1,\ldots,\phi_n-1$, denote the basis for $\FF_n$ given in Theorem
  \ref{theorem: prime power} or Theorem \ref{theorem: power of 2} and
  $f_i'$, $i=1,\ldots,\phi_{n-1}-1$, the basis for $\FF_{n-1}$. By
  Lemma \ref{lemma: projection principal to principal}, we have
  $$
    \pi_n(\div f_i)=0
  $$
  for $i=1,\ldots,\phi_n-\phi_{n-1}$, and
  \begin{equation} \label{equation: matrix}
    \begin{pmatrix}\div f_1'\\ \vdots \\
      \div f_{\phi_{n-1}-1}'\end{pmatrix}
   =\frac1{p^2}\begin{pmatrix} R & 0 \\ 0 & I\end{pmatrix}
    \begin{pmatrix}\pi_n(\div f_{\phi_n-\phi_{n-1}+1}) \\ \vdots \\
      \pi_n(\div f_{\phi_n-1}) \end{pmatrix},
  \end{equation}
  where $I$ is the identity matrix of size $\phi_{n-2}-1$ and
  $$
    R=\begin{pmatrix}1&-\beta^2&0&\cdots&\cdots&\cdots \\
    0&1&-\beta^2&\cdots&\cdots&\cdots \\
    \vdots & & & & & \vdots \\
    \cdots & \cdots & \cdots & 1 & -\beta^2 & 0 \\
    \cdots & \cdots & \cdots & 0 & 1 & -\beta^2 \\
    \cdots & \cdots & \cdots & 0 & 0 & p\end{pmatrix}
  $$
  is a square matrix of size $\phi_{n-1}-\phi_{n-2}$ whose
  superdiagonals are all $-\beta^2$ and whose diagonals are all $1$,
  except for the last one, which has $p$. Therefore, the index of
  $\pi_n(\PP_n)$ in $\PP_{n-1}$ is
  $$
    p^{2(\phi_{n-1}-1)-1}=p^{p^{n-1}(p-1)-3}.
  $$
  The structure of the factor group $\PP_{n-1}/\pi_n(\PP_n)$ can be
  easily seen from the matrix above. This completes the proof of
  the proposition.
\end{subsection}

\begin{subsection}{Proof of Proposition \ref{proposition: projection
  is an isomorphism}}
\label{subsection: projection is an isomorphism}
Consider the group homomorphism
   $$
     \pi:~\DD_n\to\pi_n(\DD_n)/\pi_n(\PP_n)=\DD_{n-1}/\pi_n(\PP_n)
   $$
    sending $D\in\DD_n$ to the coset $\pi_n(D)+\pi_n(\PP_n)$. The
    homomorphism is clearly onto, and the kernel is the group
    $\ker\pi=\PP_n+\ker\pi_n$. Thus, we have
  $$
    \DD_n/(\PP_n+\ker\pi_n)\simeq \DD_{n-1}/\pi_n(\PP_n).
  $$
  Now the group on the left-hand side is isomorphic to
  $$
    \DD_n/(\PP_n+\ker\pi_n)\simeq(\DD_n/\PP_n)\big/((\PP_n+\ker\pi_n)/\PP_n).
  $$
  Therefore, to prove that the $p$-part of $\CC_n=\DD_n/\PP_n$ is
  isomorphic to that of $\DD_{n-1}/\pi_n(\PP_n)$, it suffices to show
  that the order of $(\PP_n+\ker\pi_n)/\PP_n$ is not divisible by $p$.

  From the definition of $\pi_n$, it is easy to see that the kernel of
  $\pi_n$ is generated by divisors of the form
  $$
    D=P_{n,k}-P_{n,k+\phi_{n-1}}.
  $$
  Let $f_i$, $i=1,\ldots,\phi_n-1$, be the basis for $\FF_n$ given in
  Theorem \ref{theorem: prime power} or Theorem \ref{theorem: power of
  2}. If we write $D$ as a linear combination
  $$
    D=\sum_{i=1}^{\phi_n-1}r_i\div f_i, \qquad r_i\in\Q,
  $$
  of $\div f_i$, then the order of $D+\PP_n$ in the divisor class
  group $\CC_n$ divides the least common multiple of the
  denominators of $r_i$. We need to show that this number is not
  divisible by $p$.

  We first prove that $r_i=0$ for
  $i=\phi_n-\phi_{n-1}+1,\ldots,\phi_n-1$. By Lemma \ref{lemma:
  projection principal to principal}, we have
  \begin{equation} \label{equation: proposition 9 1}
    0=\pi_n(D)=\sum_{i=\phi_n-\phi_{n-1}+1}^{\phi_n-1}r_i
      \pi_n(\div f_i).
  \end{equation}
  Let $A=\left(\begin{smallmatrix}R&0\\0&I\end{smallmatrix}\right)$ be
  the square matrix of size $\phi_{n-1}-1$ in \eqref{equation:
  matrix}. Then we have
  $$
    \begin{pmatrix}\pi_n(\div f_{\phi_n-\phi_{n-1}+1})\\
    \vdots \\ \pi_n(\div f_{\phi_n-1})\end{pmatrix}
   =p^2A^{-1}\begin{pmatrix}\div f_1'\\ \vdots\\ \div
   f_{\phi_{n-1}-1}'\end{pmatrix},
  $$
  where $f_i'$, $i=1,\ldots,\phi_{n-1}-1$, is the basis for
  $\FF_{n-1}$ given in Theorem \ref{theorem: prime power} or Theorem
  \ref{theorem: power of 2}, and \eqref{equation: proposition 9 1} can
  be written as
  $$
   0=p^2(r_{\phi_n-\phi_{n-1}+1},\ldots,r_{\phi_n-1})A^{-1}
    \begin{pmatrix}\div f_1'\\ \vdots\\ \div f_{\phi_{n-1}-1}'
    \end{pmatrix}.
  $$
  Since $\div f_i'$ are linearly independent over $\Q$, we must have
  $$
    (r_{\phi_n-\phi_{n-1}+1},\ldots,r_{\phi_n-1})A^{-1}=(0,\ldots,0).
  $$
  It follows that $r_i=0$ for all
  $i=\phi_n-\phi_{n-1}+1,\ldots,\phi_n-1$, and
  $$
    D=\sum_{i=1}^{\phi_n-\phi_{n-1}}r_i\div f_i.
  $$

  Now, without loss of generality, we may assume that the integer $k$
  in $D=P_{n,k}-P_{n,k+\phi_{n-1}}$ satisfies $0\le
  k<\phi_n-2\phi_{n-1}$. (Let $b$ and $d$ be integers such that
  $\alpha d-bp^{n+1}=1$. Notice that if a modular unit
  $f(\tau)\in\FF_n$ has a divisor $mD$ for some integer $m$, then the
  function $f((\alpha\tau+b)/(p^{n+1}\tau+d))$ has a divisor
  $m(P_{n,k-1}-P_{n,k+\phi_{n-1}-1})$. Thus,
  $P_{n,k}-P_{n,k+\phi_{n-1}}$ and $P_{n,k-1}-P_{n,k+\phi_{n-1}-1}$
  have the same order in the divisor class group $\CC_n$.) Let $M$ be
  the square matrix of size $\phi_n-\phi_{n-1}$ whose $(i,j)$-entry is
  the order of $f_i$ at $P_{n,j-1}$. Then the order of $D$ in the
  divisor class group $\CC_n$ will divide the determinant of the
  matrix $M$. By Lemmas \ref{lemma: determinant of new part} and
  \ref{lemma: Bernoulli product} and the assumption that $p$ is a
  regular prime, the determinant of $M$ is not divisible by $p$. This
  shows that the order of $D+\PP_n$ in $\CC_n$ is not divisible by
  $p$, and therefore $|(\PP_n+\ker\pi_n)/\PP_n|$ is not divisible by
  $p$. This proves the proposition.
\end{subsection}

\begin{subsection}{Proof of Proposition \ref{proposition: kernel of
      p^2}}
\label{subsection: kernel of p^2}
  By Proposition \ref{proposition: index of projection},
  $\PP_{n-1}/\pi_n(\PP_n)$ is clearly contained in $\ker[p^2]$. Now
  suppose that $D+\pi_n(\PP_n)\in\DD_{n-1}/\pi_n(\PP_n)$ is in the
  kernel of $[p^2]$. We have $p^2D\in\pi_n(\PP_n)$. With
  \eqref{equation: p^2}, this can be written as
  $\pi_n(\iota_n(D))\in\pi_n(\PP_n)$, or equivalently
  $$
    \iota_n(D)\in\PP_n+\ker\pi_n.
  $$
  Let $f_i$, $i=1,\ldots,\phi_n-1$, be the basis for $\FF_n$ given in
  Theorem \ref{theorem: prime power} or Theorem \ref{theorem: power of
  2}. By Lemma \ref{lemma: projection
  principal to principal}, we have $\div f_i\in\ker\pi_n$ for
  $i=1,\ldots,\phi_n-\phi_{n-1}$. Hence,
  $$
    \iota_n(D)=\sum_{i=\phi_n-\phi_{n-1}+1}^{\phi_n-1}
    m_i\div f_i+D'
  $$
  for some integers $m_i$ and some divisor $D'$ in $\ker\pi_n$. Now
  notice that if we define an inner product $\gen{\cdot,\cdot}$ on
  $\DD_n$ by
  $$
    \gen{c_0P_{n,0}+c_1P_{n,1}+\cdots,d_0P_{n,0}+d_1P_{n,1}+\cdots}
   =c_0d_0+c_1d_1+\cdots,
  $$
  then for $i=\phi_n-\phi_{n-1}+1,\ldots,\phi_n-1$, $\div f_i$ is in
  the orthogonal complement of $\ker\pi_n$. The same thing is also
  true for $\iota_n(D)$ for any $D\in\DD_{n-1}$. It follows that the
  divisor $D'$ above is actually $0$ and we have $\iota_n(D)\in\PP_n$.
  Finally, by Lemma \ref{lemma: iota 2}, the fact that $\iota_n(D)$ is
  principal implies that $D$ itself is principal. This completes the
  proof of the proposition.
\end{subsection}
\end{section}

\section*{Acknowledgment}
  The authors would like to thank Professor Jing Yu for his interest
  in this work.

  Part of the work was done while the first author was visiting the
  Max-Planck-Institut f\"ur Mathematik at Bonn. He would like to thank
  the institute for providing a stimulating research environment. His
  visit was supported by a fellowship of the Max-Planck-Institut. He
  was also partially supported by Grant 96-2628-M-009-014 of the
  National Science Council, Taiwan.

  The second author was supported in part by Professor N. Yui's
  Discovery Grant from NSERC, Canada.

\end{document}